\documentclass{article}
\ifx\HCode\UnDef\else\hypersetup{tex4ht}\fi

\newcommand  \eps{ \varepsilon}

\usepackage{graphicx}
\usepackage{subfig}

\usepackage{mathrsfs}

\usepackage{amsmath}
\usepackage{amssymb}
\usepackage{textcomp}
\usepackage{enumerate}
\usepackage{hyperref}

\usepackage{comment}

\newcommand{\nn}{\nonumber}

\newtheorem{conj}{Conjecture}
\newtheorem{definition}{Definition}[section]
\newtheorem{assumption}{Assumption}[section]
\newtheorem{remark}{Remark}
\newtheorem{lemma}{Lemma}
\newtheorem{theorem}{Theorem}
\newtheorem{proposition}{Proposition}

\title{Computational singular perturbation method for nonstandard slow-fast systems}

\author{Ian Lizarraga\footnotemark[1]\
\and Martin Wechselberger\footnotemark[2]\ }
\begin{document}
\maketitle
\newcommand{\slugmaster}{%
\slugger{siads}{xxxx}{xx}{x}{x--x}}%slugger should be set to juq, siads, sifin, 

\renewcommand{\thefootnote}{\fnsymbol{footnote}}

\footnotetext[1]{School of Mathematics and Statistics, The University of Sydney, Camperdown NSW 2006, Australia (\href{mailto:ian.lizarraga@sydney.edu.au}{ian.lizarraga@sydney.edu.au}). This author acknowledges support from the ARC Discovery Project Grant DP180103022.}
\footnotetext[2]{School of Mathematics and Statistics, The University of Sydney, Camperdown NSW 2006, Australia (\href{mailto:martin.wechselberger@sydney.edu.au}{martin.wechselberger@sydney.edu.au}). This author acknowledges support from the ARC Discovery Project Grant DP180103022.}

\renewcommand{\thefootnote}{\arabic{footnote}}

\begin{abstract}
The computational singular perturbation (CSP) method is an algorithm which iteratively approximates slow manifolds and fast fibers in multiple-timescale dynamical systems. Since its inception due to Lam and Goussis \cite{lam1989}, the convergence of the CSP method has been explored in depth; however, rigorous applications have been confined to the standard framework, where the separation between `slow' and `fast' variables is made explicit in the dynamical system. This paper adapts the CSP method to {\it nonstandard} slow-fast systems having a normally hyperbolic attracting critical manifold. We give new formulas for the CSP method in this more general context, and provide the first concrete demonstrations of the method on genuinely nonstandard examples. 
\end{abstract}

%\begin{keywords} Invariant manifolds, computational methods, geometric singular perturbation theory\end{keywords}
%
%\begin{AMS} 37C10, 37M99, 37N25\end{AMS}

\pagestyle{myheadings}
\thispagestyle{plain}
\markboth{IAN LIZARRAGA AND MARTIN WECHSELBERGER}{CSP IN NONSTANDARD SLOW-FAST SYSTEMS}

\section{Introduction} \label{sec:introduction}

Most systems in nature consist of processes that evolve on disparate timescales and the observed dynamics in such systems reflect these multiple timescale features as well. Mathematical models of such multiple timescale systems are considered singular perturbation problems with slow-fast (or two timescale) problems as the most common. Models of homogeneously mixed biochemical reactions such as substrate-enzyme or ligand-receptor kinetics are prime examples.

An interesting and pervasive feature of these biochemical reaction systems is the observed transition from transient fast  kinetics to long-term slow kinetics, wherein the system settles onto a so-called {\em quasi-steady state} (QSS). Geometrically, this QSS is perceived as a lower dimensional {\em slow manifold}. Identifying such an attracting lower dimensional slow manifold provides the means to reduce the dimension of biochemical reaction systems. Such QSS reduction techniques are frequently employed in the biochemical literature with {\em Michaelis-Menten-type} laws as prime examples; see e.g.~\cite{KeenerSneyd}. 

The mathematical foundation to justify such a QSS reduction is given by Tikhonov's \cite{tikhonov1952} respectively Fenichel's \cite{fenichel1979} work on normally hyperbolic attracting slow manifolds in singular perturbation problems. Goeke, Noethen and Walcher \cite{nowa,goeke2014} provide a comprehensive discussion on the general setup of Fenichel's geometric singular perturbation theory (GSPT) with an emphasis on explaining when a QSS reduction is justified or when it leads to erroneous results.

Schauer and Heinrich \cite{schauer1983} explored a homotopy approach to find QSS approximations in biochemical reaction networks, by formulating a perturbation problem in terms of a small parameter $\eps>0$ amplifying the ratio of magnitudes of slow and fast reaction rates.\footnote{We note that a small perturbation parameter $\eps>0$ could be properly identified via dimensional analysis.} After using stoichiometry to derive the dynamical system from the network, the fast reactions are grouped into a vector $W$, and slow into a vector $V$:
 \begin{eqnarray}
z' &=& NW(z) + \eps RV(z), \label{eq:stiefsplitting}
\end{eqnarray}
with state vector $z\in\mathbb{R}^n$, $N$ and $R$ are constant stoichiometric matrices with full column rank $m_N,m_R<n$, and $\eps>0$ is the homotopy parameter.  
With this splitting, a  $k$-dimensional manifold of stationary states (QSS), $1\le k<n$, can be identified in the singular limit $\eps\to 0$,
\begin{eqnarray}
S = \{z\in\mathbb{R}^n\;:\; NW(z) = 0 \Rightarrow W(z) = 0\}\,.
\end{eqnarray}
We emphasize here that the slow-fast splitting \eqref{eq:stiefsplitting} is {\it nonstandard} from the point of view of GSPT %multiple-timescale dynamical systems theory, 
in the sense that slow-fast reactions are distinguished rather than slow-fast variables. Nevertheless, the theory of Fenichel \cite{fenichel1979} still applies: under the appropriate geometric condition (see Section~\ref{sec:nongspt} for details), an invariant slow manifold $S_\eps$ perturbs from $S$ for sufficiently small $\eps > 0$. 
%This splitting arises in other contexts, including stick-slip friction models \cite{hinrichs1998}. 
%
\begin{figure}[t]
  \centering
        (a)  \includegraphics[height=0.35\textwidth,width=0.95\textwidth]{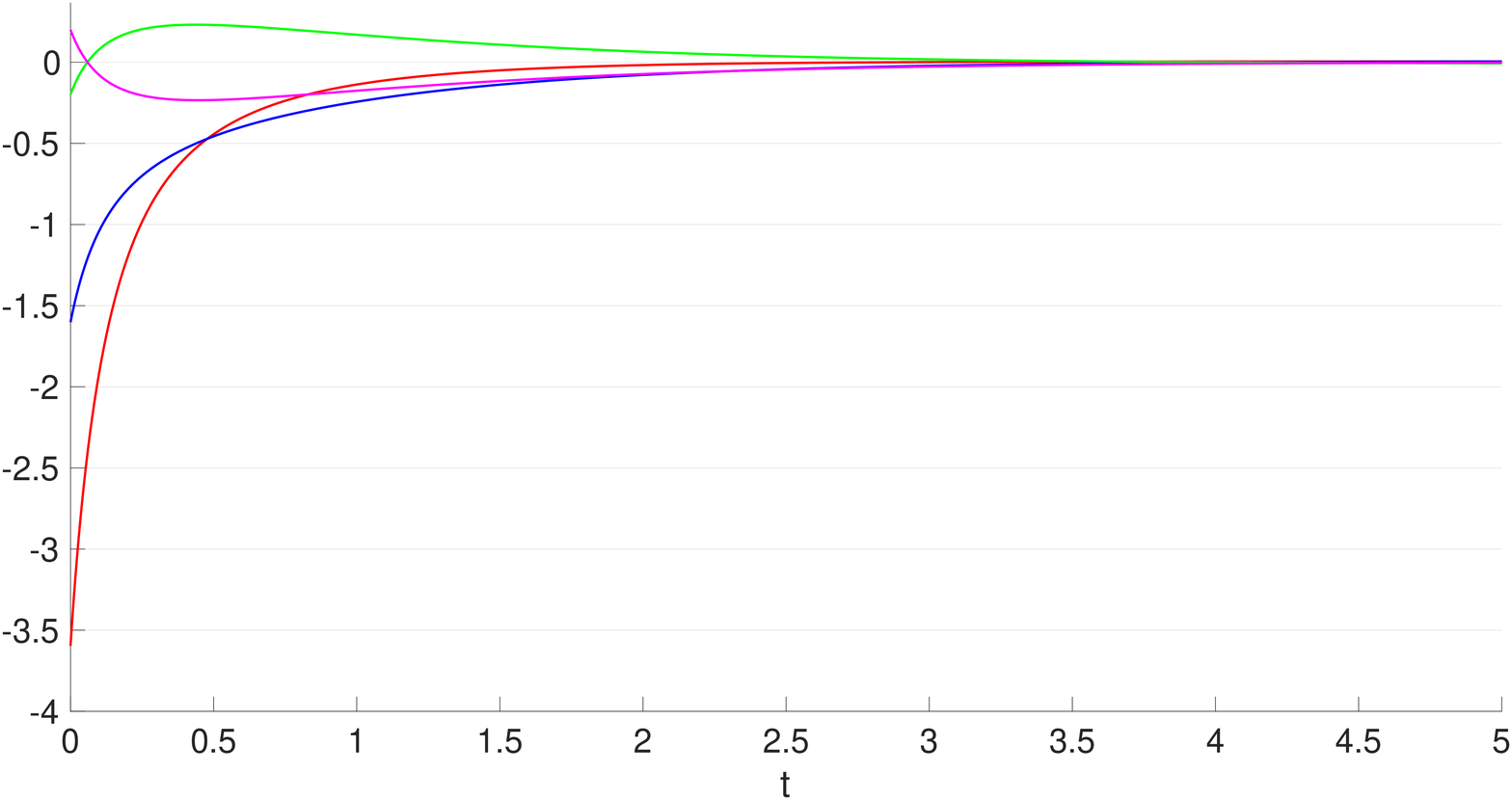}\\
(b)    \includegraphics[width=0.95\textwidth]{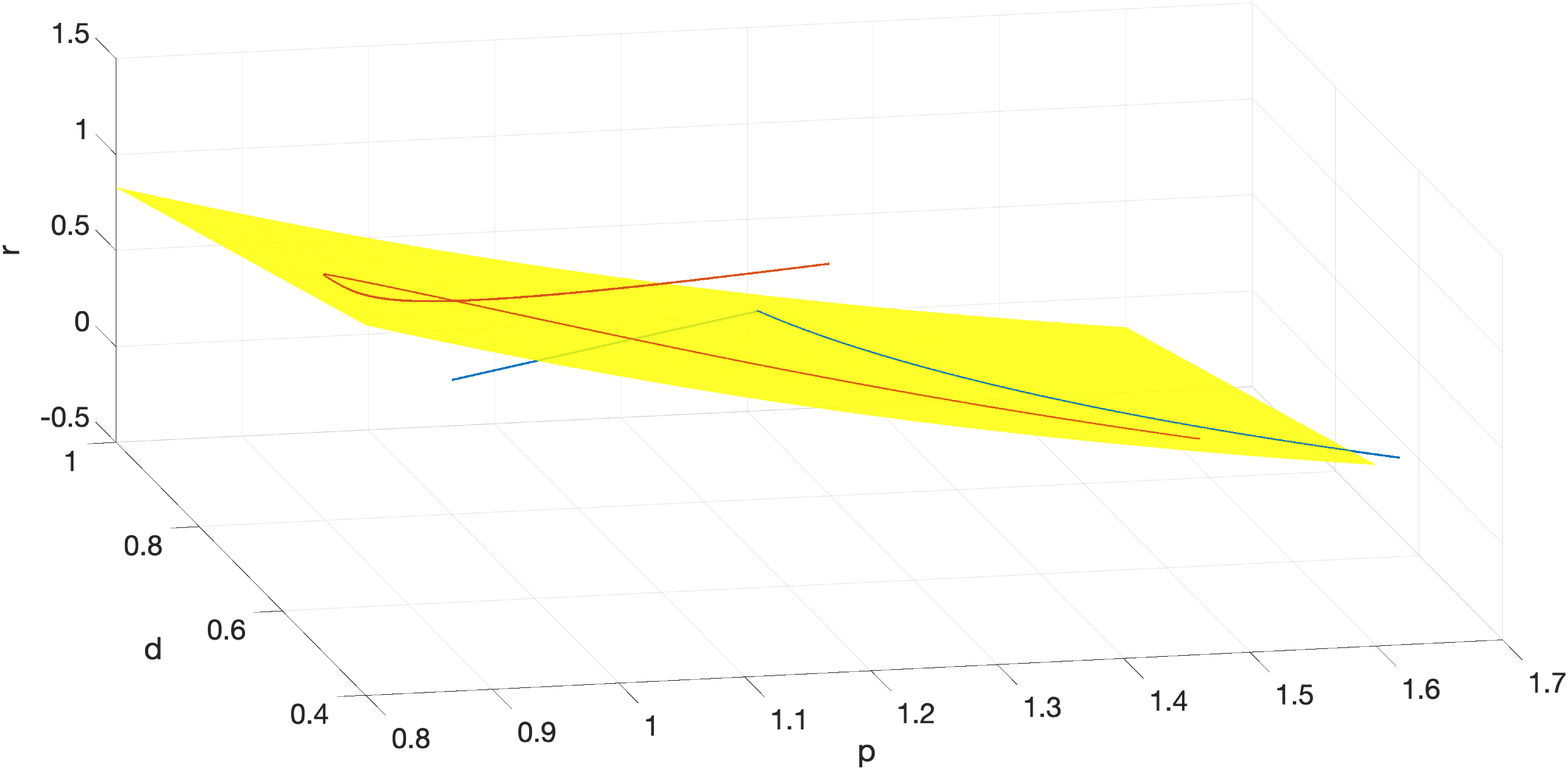}
      \caption{System \eqref{eq:slime2} with $\eps = 0.01$, $k_i = S = 1$, $c = 3$ and initial conditions: $(p,d,r,b) = (1,0.8,0.6,0.4)$;
      (a) Time series of the velocities $p',d',r',b'$ (respectively in red, blue, green, magenta). (b) Projection onto $(p,d,r)$ space of two trajectory segments together with a portion of the manifold of equilibria defined by \eqref{eq:slimemanifold}. Initial conditions: blue, $(p,d,r,b) = (1,0.8,0.6,0.4)$; red, $(p,d,r,b) = (1.2,0.5,1.4,0.4)$. Trajectories computed using a Dormand-Prince ODE solver in MATLAB R2018a \cite{matlab}.} 
 
      \label{fig:slime}
\end{figure}

As a motivating example, consider the following four-dimensional slime mold cell communication model \cite{goeke2013}:

\begin{eqnarray}
\begin{pmatrix}
p' \\ d' \\ r' \\ b'
\end{pmatrix} &=& \begin{pmatrix}
2 & 2 \\
0 & 1 \\
1 & 0 \\
-1 & 0
\end{pmatrix} \begin{pmatrix}
-k_5 rp^2 + k_{-5}b\\
-k_4 dp^2 + k_{-4}(c-d-r-b)
\end{pmatrix} + \eps \begin{pmatrix}
k_3 - k_{-3}p + k_2 Sb \\
-k_1 d + k_{-1}r \\
k_1 d - k_{-1} r\\
0
\end{pmatrix}. \label{eq:slime2}
\end{eqnarray}
Here, $p$ refers to the concentration of cAMP; $d$ and $r$ represent transmembrane receptors; and $b$ refers to the bound state of $r$. The parameters $k_i \geq 0$ are constant reaction rates, the parameter $c\geq 0$ represents a conserved quantity that arises from the model reduction of the original five-dimensional model due to Stiefenhofer \cite{stiefenhofer}, and the parameter $S \geq 0$ denotes the concentration of ATP which is assumed to be constant. This system contains a two-dimensional manifold of equilibria given as a graph
\begin{equation} \label{eq:slimemanifold}
S = \{(p,d,r,b)\in\mathbb{R}^4: (r,b) = \psi(p,d)\} 
\end{equation}
with
\begin{equation}
\psi(p,d) = 
\left( \frac{k_{-5}(-k_4 dp^2 + k_{-4} (c-d))}{k_4 ( k_5 p^2 + k_{-5})}\, ,\,\frac{k_5 p^2 (-k_4 dp^2 + k_{-4} (c-d))}{k_4 ( k_5 p^2 + k_{-5})}\right)\,. 
\end{equation}
When $\eps>0$ is sufficiently small, numerical observations reveal the decay of fast transient motion as trajectories are attracted to a low-dimensional invariant set $S_\eps$ near $S$, as shown in Figure \ref{fig:slime}. Our objective is to approximate such invariant sets arbitrarily well.

Several algorithmic \& computational techniques have been developed to identify these QSS (i.e.~a lower dimensional slow manifold),  and to identify the simplified system underlying long-term stability.  
Lam and Goussis  \cite{lam1989,lam1994} devised the {\it computational singular perturbation} (CSP) method for chemical kinetics modelling which is an iterative procedure that generates a sequence of CSP manifolds and CSP fibers which approximate the slow manifolds and fast fibers of the dynamical system; see Section~\ref{sec:csp2} for details. The convergence of these CSP objects has been analyzed by Kaper, Kaper, and Zagaris in a series of papers \cite{kkz2004a,kkz2004b,kkz2005,kkz2015}. They rigorously prove this convergence in the case of slow-fast vector fields satisfying a spectral gap condition for the eigenvalues of the Jacobian of the vector field along an attracting invariant manifold, but focus on examples in the standard form. Indeed, the practical applicability of the CSP method has remained unclear in the case where we cannot easily identify the slow or fast components of the vector field {\it a priori}. Lam and Goussis \cite{lam1994} give an example of CSP applied to a linearized problem by populating the initial conditions with (constant) eigenvectors, but pointed out the necessity of a modification for more difficult nonlinear problems where the Jacobian of the vector field is point-dependent.
 
Our approach is to avoid this defect for nonstandard slow-fast systems of the form
\begin{equation}\label{eq:slow-fast-general}
z' = H(z,\eps)= N(z)f(z)+\eps G(z,\eps)  
\end{equation}
which includes system \eqref{eq:slime2}. % as well as standard slow-fast systems. 
The idea is that the singular limit vector field factorization of \eqref{eq:slow-fast-general}, i.e.~the vector field $H(z,0)=N(z)f(z)$, encodes a wealth of geometric information. We will make explicit use of this factorization to provide formulas for a nonstandard version of the CSP method.  
% (REMARK: mentioned that somewhere else?) We remark that there are known slow-fast models for which such leading-order factorizations have not been given, such as the Olsen model of oxidase-peroxidase reaction. One representative analysis of the Olsen model from the point of view of blow-up theory is given in \cite{kuehnolsen2015}.

Although this paper is concerned with the CSP method, there are in fact a variety of iterative schemes to approximate invariant manifolds; the zero-derivative principle (ZDP), intrinsic low-dimensional manifolds (ILDM), and center manifold normal-form reductions \cite{roberts2008} are among a few of these (chapter 11 in \cite{kuehn2015}  give an overview of the first three methods). Many of these schemes are interconnected; for instance, the CSP and ZDP iterations differ only by a multiplicative factor \cite{kkz2005}. 

The paper is organised as follows: In Sec.~\ref{sec:nongspt}, we describe a general framework for nonstandard slow-fast systems. In Sec.~\ref{sec:csp2} we define the two-step CSP update. Our main results are provided in Secs.~\ref{sec:csp3} and \ref{sec:eg}: we give specific formulas for initializing the CSP method in the nonstandard context, study the output of the first update carefully, and then give examples demonstrating these formulas. We conclude in Sec.~\ref{sec:conc} by highlighting some fruitful new connections between the CSP method and the factorization given in the second section. We also provide an appendix (App.~\ref{sec:csp1}), which expands on the framework underlying the CSP update step.

%\newpage
\section{Nonstandard slow-fast dynamical systems}\label{sec:nongspt}

We begin by giving an abbreviated treatment of a general framework for nonstandard slow-fast systems. Much of this material in fact appears in Fenichel's seminal work on GSPT \cite{fenichel1979}; see also \cite{nowa,goeke2014}. This approach has been further developed by Wechselberger \cite{wechselberger2019} and extends the framework to loss of normal hyperbolicity. 

We are interested in two-timescale (or {\it slow-fast}) dynamical systems of the form \eqref{eq:slow-fast-general}, which we restate here for convenience:
\begin{equation} \label{eq:master}
z' = \frac{dz}{dt} = H(z,\eps) = N(z)f(z)+\eps G(z,\eps)\,,
\end{equation}
with state variable $z \in \mathbb{R}^n$, $n\times (n-k)$ matrix $N(z)$ formed by column vectors $N^i(z)=(N^i_1(z),\ldots,N^i_n(z))^\top$ with sufficiently smooth functions $N^i:\mathbb{R}^{n}\to\mathbb{R}$, $i=1,\ldots, n-k$, $f(z)=(f_1(z),\ldots,f_{n-k}(z))^\top$ a column vector of sufficiently smooth functions $f_i:\mathbb{R}^{n}\to\mathbb{R}$, $i=1,\ldots,n-k$, $G(z,\eps)=(G_1(z,\eps),\ldots,G_{n}(z,\eps))^\top$ a column vector of sufficiently smooth functions $G_i:\mathbb{R}^{n}\to\mathbb{R}$, $i=1,\ldots,n$, and $\eps \ll 1$ characterizes the ratio of timescales in the system. \\

\begin{definition}
~Let $S_0$ denote the set of equilibria of system \eqref{eq:master} in the singular limit $\eps \to 0$.
If there exists a subset $S\subseteq S_0$ which forms a $k$-dimensional differentiable manifold of equilibria with $1\le k< n$, then system \eqref{eq:master} defines a {\em singular perturbation problem}. \\
\end{definition}

\begin{assumption}\label{ass1}
~ System \eqref{eq:master} is a singular perturbation problem with a single subset $S\subseteq S_0$,
\begin{equation}\label{eq:critman}
S = \{z\in \mathbb{R}^n: f(z) =  0\}\,, 
\end{equation}
which forms  a $k$-dimensional differentiable manifold of equilibria, $1\le k< n$, called the  {\rm critical manifold}.\\
\end{assumption}

\begin{assumption}\label{assN}
~In system \eqref{eq:master}, the matrix $N(z)$ has full (column) rank for all $z\in S$.\\
\end{assumption}

\noindent
Next consider system \eqref{eq:master}, rescaled from the fast timescale $t$ to the slow timescale $\tau = \eps t$:
\begin{equation} \label{eq:slowsystem}
\dot{z} = \frac{dz}{d\tau} = \frac{1}{\eps}H(z,\eps) = \frac{1}{\eps}N(z) f(z) +  G(z,\eps),
\end{equation}
Systems \eqref{eq:master} and \eqref{eq:slowsystem} are equivalent when $\eps > 0$ but their singular limits $\eps \to 0$ are not. In fact, they carry complementary, lower dimensional information, and it is a cornerstone of GSPT to {\em concatenate} the information from these two limiting problems to deduce the dynamics of the full system  \eqref{eq:master} respectively \eqref{eq:slowsystem}.\\

\subsection{The layer problem} \label{layerproblem}

%A cornerstone of  geometric singular perturbation theory (GSPT) analysis is Fenichel's theorem. The standard approach is to study  concatenations of trajectory segments obtained from solutions of reduced subsystems, the {\it layer} and {\it reduced} problems, when $\eps = 0$. 
We focus first on system \eqref{eq:master} evolving on the fast time scale $t$.\\

\begin{definition} 
~The {\rm layer problem} of system \eqref{eq:master} is given by the formal limit $\eps \to 0$:
\begin{equation}	
z' = H(z,0) = N(z)f(z).
\end{equation}
\end{definition}

\noindent
Under Assumption~\ref{ass1}, the set $S$ forms a $k$-dimensional manifold of equilibria of the layer problem. Hence, the Jacobian $Dh|_S$ along $S$ has $k$ {\em trivial} zero eigenvalues and $(n-k)$ {\em nontrivial} eigenvalues.\\
%since $S$ forms a manifold of equilibria (of dimension $k$) for the layer problem.

\begin{definition}
~A $k$-dimensional critical manifold $S$ is called {\rm normally hyperbolic} if the $(n-k)$ nontrivial eigenvalues of the Jacobian $Dh|_{S}=(N D\!f)|_{S}$ are bounded away from the imaginary axis. \\
\end{definition}

\noindent

The existence of a normally hyperbolic manifold implies the following splitting:
\begin{eqnarray}
T_z\mathbb{R}^n &=& T_z S \oplus \mathcal{N}_z,\quad\forall z\in S \label{eq:layer}
\end{eqnarray}
where $T_z S$ denotes the tangent space of the critical manifold $S$ at $z$, coinciding with the kernel of the linear map $DH(z,0)$, and $\mathcal{N}_z$ is the unique complement of the splitting identified with the quotient space $T_z \mathbb{R}^n / T_z S$.  We call $\mathcal{N}_z$ the {\it linear fast fiber with basepoint $z \in S$}. Repeating this construction across all points of $S$, we obtain the tangent bundle $TS=\cup_{z\in S} T_zS$ of $S$, and the transverse linear fast fiber bundle $\mathcal{N}=\cup_{z\in S} \mathcal{N}_z$. This splitting induces unique projection operators onto the tangent bundle $TS$ and the linear fast fibre bundle $\mathcal{N}$ as follows:

\begin{eqnarray}
\Pi^S &:& TS \oplus \mathcal{N} \to TS. \label{eq:projS}\\
\Pi^N &:& TS \oplus \mathcal{N} \to \mathcal{N}  = I - \Pi^S. \label{eq:projN}
\end{eqnarray}

\begin{lemma} \label{lemma2}
~At any point $z \in S$ of a normally hyperbolic manifold $S$, the rows of $Df(z)$ form a basis of $T_zS^{\bot}$ and the columns of $N(z)$ form a basis of $\mathcal{N}_z$.\\
\end{lemma}

{\it Proof.} 
The manifold $S$ is the zero level set of $f(z)$ which forms a $k$-dimensional manifold, i.e. $D\!f(z)|_{S}$ has full row rank $(n-k)$ for any $z\in S$. Thus the rows of $Df(z)$ form a basis of the orthogonal complement of $T_zS$. 

Since $S$ is normally hyperbolic, we have $DH(z,0)$ has rank $(n-k)$ which implies $\ker DH(z,0) =  T_zS$ at points $z \in S$. Under Assumption~\ref{assN}, $N(z)$ has full column rank $(n-k)$. Thus the column spaces of $DH(z,0)$ and $N(z)$ coincide at points $z\in S$ and, hence, the columns of $N(z)$ form a basis of $\mathcal{N}_z$.  \hfill $\Box$\\

\begin{comment}
By construction, the column space of the Jacobian forms a basis of the fibers with basepoints on $S$. The Jacobian $DH(z,0)$ has rank $(n-k)$ at points $z \in S$, and thus the column spaces of $DH(z,0)$ and $N(z)$ coincide by the previous lemma. 

The subspace ker $DH(z,0)$ has dimension $k$. Let $W(z) = [W_1(z)~W_2(z)~\cdots~W_k(z)]$ denote the $n\times k$ matrix formed by concatenating a basis set of ker $Dh$. Thus,

\begin{eqnarray*}
DH(z,0) W(z) &=& 0\\
N(z) Df(z) W(z) &=& 0\\
\Rightarrow Df(z) W(z)&=&0,
\end{eqnarray*}

where the final implication comes from $N(z)$ having full column rank.  \hfill $\Box$\\
\end{comment}

\begin{lemma} \label{lemma3}
~Let $S$ be a normally hyperbolic manifold. Then the $(n-k) \times (n-k)$ square matrix $D\!f N|_{S}$ is regular, and its eigenvalues are equal to the set of nontrivial eigenvalues of $Dh|_{S}$. \\
\end{lemma}

{\it Proof.} Let the linear operator $DH(z,0)|_{S}$ act on the basis of the invariant subset $\mathcal{N}$, i.e.
$$
DH(z,0)N(z)|_S = (N(z) D\!f(z)) N(z)|_S = N(z) (D\!f(z) N(z))|_S.
$$
By Lemma~\ref{lemma2}, the $(n-k) \times (n-k)$ square matrix $D\!f N|_{S}$ is necessarily a regular matrix with $(n-k)$ nonzero eigenvalues which coincide with the nontrivial eigenvalues of $DH(z,0)|_S$. Note, if $p(z)\in\mathbb{R}^{n-k}$ is an eigenvector of  $D\!f N|_{S}$ with eigenvalue $\mu(z)$ then $N(z)p(z)\in\mathbb{R}^{n}$ is the corresponding eigenvector of $DH(z,0)$ with the same eigenvalue.
\hfill $\Box$\\

\begin{assumption} \label{assu1b} 
%{(Normally hyperbolic attracting critical manifold)}\\
 ~The $k$-dimensional critical manifold $S$ of system \eqref{eq:master} is normally hyperbolic and attracting, i.e. all $(n-k)$ nontrivial eigenvalues have negative real part.\\ 
\end{assumption}

\begin{remark}
~System \eqref{eq:slime2} satisfies Assumptions~\ref{ass1}--\ref{assu1b}. A rich source of examples comes from chemical reaction networks; Schauer and Heinrich \cite{schauer1983} provide a derivation of such models with the goal of identifying quasi-steady state approximations in biochemical reaction networks.  

Nevertheless, it is unclear whether general techniques exist to compute such factorizations in a typical slow-fast system. In the case of polynomial vector fields, Goeke and Walcher \cite{goeke2014} have used division algorithms for algebraic varieties to factor $H(z,0)$.\\
\end{remark}

\begin{remark} \label{remark:isolated}
~In this paper we do not concern ourselves with singularities of $H(z,0)$ outside of the critical manifold $S$; we only require the local geometric structure provided by the differentiable manifold and transverse fibers in the proceeding arguments. A simple example of a system having an isolated singularity as well as a critical manifold is given in Sec. \ref{sec:circle}.\\
\end{remark}

\subsection{The reduced problem} \label{reducedproblem}

Now consider system \eqref{eq:slowsystem} evolving on the slow timescale $\tau = \eps t$.
%\begin{eqnarray}
%\dot{z} &=& \frac{1}{\eps}N(z) f(z) +  G(z,\eps), %\label{eq:slowsystem}
%\end{eqnarray}
%where $\dot{ } = d/d\tau$. The equations \eqref{eq:master} and \eqref{eq:slowsystem} are equivalent when $\eps > 0$; however,  
The singular limit $\eps \to 0$ of system \eqref{eq:slowsystem} requires more care, i.e.~this limit will only be well-defined provided that we restrict the phase space to $S$ and that $G(z,0)$ is restricted to the tangent bundle $TS$. The projection operator $\Pi^S$ \eqref{eq:projS} allows us to formulate this limit.\\

\begin{definition}
~The {\rm reduced problem} of system \eqref{eq:slowsystem} is
\begin{eqnarray}
\dot{z} &=& \frac{d}{d\tau} z = \Pi^S \left.\frac{\partial}{\partial \eps} H(z,\eps)\right|_{\eps = 0} = \Pi^S G(z,0). \label{eq:reduced}
\end{eqnarray}
\end{definition} 
%
%These projections are not necessary to formulate the CSP method; however, we highlight in Sec. \ref{sec:proj} an intriguing connection between the algorithm and oblique projection operators defined by the vector field. \\
% 

%We end this section by using these lemmas to give a well-defined formula for the projector \eqref{eq:projS}, which is used to define the reduced system.\\

\begin{figure}[!]
  \centering
    \includegraphics[width=0.8\textwidth,height=0.4\textheight]{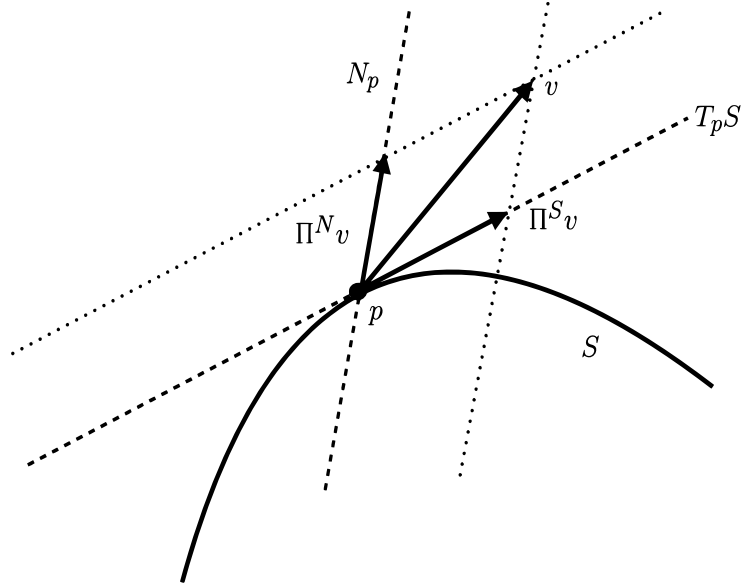}
      \caption{A sketch of the projectors $\Pi^S$ \eqref{eq:projS} and $\Pi^N$ \eqref{eq:projN} defined at a point $p \in S$. The invariant subspaces $T_pS$ and $N_p$ are illustrated by dashed lines, together with an arbitrary vector $v \in T_p\mathbb{R}^n$. The dotted lines are parallel translates of these subspaces along the oblique projections $\Pi^S v$ and $\Pi^N v$.}
      \label{fig:projection1}
\end{figure}

\begin{definition} \label{def:proj}
Consider the splitting $\mathbb{R}^n = \mathcal{V} \oplus \mathcal{W}$, where $\mathcal{V}$ has dimension $n-k$ and $\mathcal{W}$ has dimension $k$. The {\rm oblique projection of a vector $p = x+y$, $x \in \mathcal{V}$ and $y \in \mathcal{W}$, onto the subspace $\mathcal{V}$ parallel to $\mathcal{W}$} is a linear map $\Pi^{\mathcal{V}} = (\Pi^{\mathcal{V}})^2$ satisfying $\Pi^{\mathcal{V}}(p)=x$.
%\begin{itemize}
%\item For each $x \in \mathcal{V}$, we have $\Pi^{\mathcal{V}}(x) = x$,
%\item For each $y \in \mathcal{W}$, we have $\Pi^{\mathcal{V}}(y) = 0$,
%\item For each $p \in \mathbb{R}^n$, we have $\Pi^{\mathcal{V}}(p) \in \mathcal{V}$.\\
%\end{itemize}

Suppose $V$ (resp. $U$) is an $n \times (n-k)$ matrix whose column vectors span $\mathcal{V}$ (resp. $\mathcal{W}^{\bot}$). Then it can be shown that $\Pi^{\mathcal{V}}$ has the following matrix representation:
\begin{eqnarray}
\Pi^{\mathcal{V}} &=& V(U^{\top} V)^{-1}U^{\top}. \label{eq:obliqueproj}
\end{eqnarray}
The (complementary) {\rm oblique projection onto the subspace $\mathcal{W}$ parallel to $\mathcal{V}$} is given by
\begin{eqnarray}
\Pi^{\mathcal{W}} &=& I - \Pi^{\mathcal{V}}. \label{eq:obliqueprojcomp}
\end{eqnarray}
\end{definition}

\noindent
In the present context, the splitting $T\mathbb{R}^n = \mathcal{N} \oplus TS$ induces an oblique projection onto the tangent bundle of the critical manifold, parallel to the fast fibers (see Fig. \ref{fig:projection1}). Lemma \ref{lemma2} gives us a matrix representation using the matrices $N$ and $Df$:
\begin{eqnarray}
\Pi^S &=& I - N(D\!fN)^{-1}Df.\label{eq:projSformula}
\end{eqnarray}
By Lemma \ref{lemma3}, the matrix $D\!fN$ is regular and, hence, the inverse is well defined. Equivalent projection formulas appear in \cite{fenichel1979} and \cite{goeke2014}.\\

%The following assumption will help us to organize the rich geometric structure of the system when $\eps = 0$. \\

\begin{comment}
%
\begin{assumption}{(Factorization of $H(z,0)$) } \label{assu2}\\
We assume that the leading order term of \eqref{eq:master} may be factorized on $\mathbb{R}^n$ as 
\begin{eqnarray}
H(z,0) &=& N(z) f(z),
\end{eqnarray}

where $N(z)$ is an $n \times (n-k)$ matrix of full column rank, formed by columns of sufficiently smooth functions, and $f(z)$ is $k$-dimensional column vector of sufficiently smooth functions. We further demand that the critical manifold $S$ be defined as the following regular level set:

\begin{eqnarray}
S &=& \{z \in \mathbb{R}^n: f(z) = 0\}. \label{eq:fcrit}
\end{eqnarray}

This implies that $Df(z)$ has full row rank. Finally, we assume that singularities of $H(z,0) = N(z) f(z)$ are isolated on $z\notin S$, if they exist.  
\end{assumption}
~\\

~\\Per Assumptions \ref{ass1}--\ref{assu1b}, we focus on dynamical systems which can be written in the form

\begin{eqnarray}
z' &=& N(z) f(z) + \eps G(z,\eps).  \label{eq:master2}
\end{eqnarray}

\end{comment}

\begin{remark}
%We also briefly comment on  {\it standard} systems. 
~It is easily shown that {\em standard} singular perturbation problems are a special case of a nonstandard problem \eqref{eq:master}. Indeed, select
 %the factorization given in Assumption \ref{assu2}; 
%
\begin{eqnarray}
N &=& \begin{pmatrix}
\mathbb{O}_{k,n-k}\\ \mathbb{I}_{n-k,n-k}
\end{pmatrix} \label{eq:trivial}
\end{eqnarray}
and describe the variables and vector field in terms of components  $z = (x,y) \in \mathbb{R}^{k} \times \mathbb{R}^{(n-k)}$ and $G(x,y,\eps) = (g(x,y,\eps),\tilde{f}(x,y,\eps))$. Then $z' = N(z)f(z) + \eps G(z,\eps)$ gives
\begin{equation}\label{eq:standard}
\begin{aligned}
x' &= \eps g(x,y,\eps)\\
y' &= f(x,y,\eps) + \eps \tilde{f}(x,y,\eps)\,, 
\end{aligned}
\end{equation}
which is a standard singular perturbation problem.
From \eqref{eq:standard} we observe that $x$ lists the $k$ slow variables and $y$ lists the $(n-k)$ fast variables. These systems lie in contrast to the {\it nonstandard} case, where one cannot expect such a trivial factorization \eqref{eq:trivial} to exist globally.
\end{remark}

\subsection{Fenichel Theory} \label{sec:fenichel}

In the case of normally hyperbolic critical manifolds, QSS reductions onto a slow invariant manifold are justified by the following theorem. \\ %In particular, we have

\begin{theorem}[Fenichel's Theorem \cite{fenichel1979,jones1995,kuehn2015}] \label{thm:fenichel}
~Given the system \eqref{eq:master} with a $C^r$-smooth vector field and a compact normally hyperbolic critical manifold, the following hold for $\eps > 0$ sufficiently small: 
\begin{itemize}
\item There exists a locally invariant $C^r$-smooth, normally hyperbolic slow manifold $S_{\eps}$ that is $C^r$ $\mathcal{O}(\eps)$-close to $S$.
\item The flow on $S_{\eps}$ converges to the reduced flow on $S$ as $\eps \to 0$.
\item There are $C^r$-smooth locally invariant stable and unstable manifolds, $\mathcal{F}^s(S_{\eps})$ and $\mathcal{F}^u(S_{\eps})$.
\item These manifolds admit nonlinear, $C^{r-1}$-smooth foliations $\{\mathcal{F}^s(p): p \in S_{\eps}\}$ resp. $\{\mathcal{F}^u(p): p \in S_{\eps}\}$. Furthermore, these families are positively (resp. negatively) invariant on fibers; i.e. if $\phi_t$ is the time-$t\geq 0$ flow map of \eqref{eq:master} and if $p, \phi_t(p) \in S_{\eps}$, then
\begin{eqnarray*}
\phi_t(\mathcal{F}^s(p)) &\subseteq& \mathcal{F}^s(\phi_t(p)).
\end{eqnarray*}
An analogous statement holds for the fibers $\mathcal{F}^u(p)$, with $p \in S_{\eps}$. 
\item There exist constants $C, \lambda > 0$ such that if $q \in \mathcal{F}^s(p)$, then for $t \geq 0$ we have
\begin{eqnarray*}
|| \phi_{t}(p) - \phi_{t}(q) || &<& C e^{-\lambda t}.
\end{eqnarray*}
Analogous rate estimates hold for the fibers $\mathcal{F}^u(p)$, with $p \in S_{\eps}$. 
\end{itemize}
\end{theorem}

The invariant slow manifold $S_{\eps}$ and its invariant nonlinear foliation $\mathcal{F}_{\eps}$ (which split into the sub-foliations given in the theorem above) organise the dynamics of the full system: trajectories perturb from concatenated orbits of the layer and reduced problems under the appropriate transversality conditions. In analogy to the $\eps = 0$ setting, we denote the linear fast fiber at basepoint $p \in S_{\eps}$  by $\mathcal{N}_{\eps,p}$ and the linear fast fiber bundle by $\mathcal{N}_{\eps} = \cup_{p \in S_{\eps}} \mathcal{N}_{\eps,p}$.

%\newpage
\section{The CSP iteration}\label{sec:csp2}

Suppose that the invariant slow manifold $\mathcal{S}_{\eps}$ of a nonstandard slow-fast system \eqref{eq:master} is (locally) given by the graph of a function $y = h_{\eps}(x)$ with $h_{\eps}: \mathbb{R}^{k} \to \mathbb{R}^{n-k}$. We express $h_{\eps}(x)$ as an asymptotic series in $\eps$:
\begin{eqnarray}
y &=& h_{\eps}(x) = h_0(x) + \eps h_1(x) + \eps^2 h_2(x) + \cdots + \eps^{j}h_j(x) + \mathcal{O}(\eps^{j+1}).
\end{eqnarray}
The terms $h_{j}(x)$ may be obtained reinserting this series into \eqref{eq:master} and matching coefficients in orders of $\eps$. The CSP method adopts a different, iterative approach to efficiently compute not only the terms $h_{j}(x)$, but also a similar asymptotic approximation to the linear fast fibers transverse to the slow manifold. The motivating idea is to understand how the variational equation of a dynamical system is affected by changes in basis. 
%
\begin{comment}
There is a clever choice of basis which explicitly characterizes the invariant slow manifold and fast fibers, and the CSP method iteratively approximates this special basis.

In this section, we first briefly review the geometric framework underlying the CSP method. Then we define the CSP objects which approximate the invariant objects in the dynamical system. Finally, we define the one-step and two-step CSP methods explicitly. Auxiliary lemmas and convergence theorems are relegated to the appendix where appropriate. 
\end{comment}

\subsection{Geometric framework of the CSP method} Consider a smooth vector field $H(z)$, to which we append the variational equation to produce the dynamical system on the tangent bundle:
\begin{eqnarray}
z' &=& H(z)\\
H' &=& DH(z) H(z).\nonumber
\end{eqnarray}
Under some arbitrarily chosen splitting $z = (x,y)$, where $x \in \mathbb{R}^{n-k}$ and $y \in \mathbb{R}^{k}$, we can write the variational equation in block component form:
\begin{eqnarray}
\begin{pmatrix}
H_x' \\ H_y'
\end{pmatrix} &=& 
\begin{pmatrix}
D_{x} H_x & D_{y} H_x \\
D_{x} H_y & D_{y} H_y
\end{pmatrix} \begin{pmatrix}
H_x \\ H_y
\end{pmatrix}. \label{variational}
\end{eqnarray}

Now suppose $A(z) = [A_f(z)~A_s(z)]$ is a smooth, regular  $n\times n$ matrix for all $z\in S$ with $B(z) = \begin{pmatrix} B_{s\bot}(z) \\ B_{f\bot}(z)\end{pmatrix}$  its dual. Here, the first and second block columns of $A$, denoted $A_f$ and $A_s$, are of sizes $n \times (n-k)$ and $n \times k$, respectively; similarly, the first and second block rows $B_{s\bot}$ and $B_{f\bot}$ have sizes $(n-k) \times n$ and $k \times n$, respectively. The vector field $H(z)$ expressed in this new basis is given by
\begin{eqnarray}
g(z) &=& B(z)H(z).\label{eq:newbasis}
\end{eqnarray}
\begin{definition}
~For a smooth $n\times m$ matrix function $X$ and a smooth vector field $H$, $[X,H]$ is denoted the {\rm Lie bracket of $X$ and $H$} and defined as the $n\times m$ matrix whose $i$th column is
\begin{eqnarray*}
[X,H]_i &=& [X_i,H]\\
&=& (DH)X_i - (DX_i)H,
\end{eqnarray*}
where $i = 1, \cdots, m$ and $A_i$ refers to the $i$th column of $A$.\\
\end{definition}

It can be shown that the variational equation of the transformed vector field can be written compactly in the form
\begin{eqnarray}
g' &=& \Lambda(A,B,H) g\,, \label{eq:lambdaeq}
\end{eqnarray}
where the nonlinear operator $\Lambda(A,B,H)$ has the algebraic structure of a Lie bracket:
\begin{eqnarray}
\Lambda(A,B,H) &=& B[A,H]\,; \label{eq:liebracket}
\end{eqnarray}
see the appendix (Sec. \ref{app:lie}) for a derivation of this result. The operator $\Lambda$ can be written in block-component form:
\begin{eqnarray}
\Lambda &=& \begin{pmatrix}
\Lambda_{ff} & \Lambda_{fs}\\
\Lambda_{sf} & \Lambda_{ss}
\end{pmatrix} = \begin{pmatrix}
B_{s\bot}[A_f,H] & B_{s\bot}[A_s,H]\\
B_{f\bot}[A_f,H] & B_{f\bot}[A_s,H]
\end{pmatrix}. \label{eq:lambda}
\end{eqnarray}
The key insight is that in the presence of a $k$-dimensional invariant manifold $\mathcal{M}$ and corresponding transverse linear fiber bundle $\mathcal{N}$, the vector field $H$ may be expressed in a clever choice of basis $A(z)$ and dual $B(z)$ which block-diagonalizes $\Lambda$. As a consequence, the  invariant manifold $\mathcal{M}$ and the linear fiber bundle $\mathcal{N}$ are easily characterized in terms of this basis as follows:
\begin{eqnarray}
\mathcal{M} &=& \{z\in \mathbb{R}^n: B_{s\bot}(z)H(z) = 0\}, \label{eq:manifold} \\
\mathcal{N} &=& \bigcup_{p\in \mathcal{M}}\mathcal{N}_p = \bigcup_{p\in \mathcal{M}}\text{Col}(A_f(p))\,; \label{eq:fiber}
\end{eqnarray}
see the appendix (Sec. \ref{app:lambda}) for details. This formalism applies directly to the case of an invariant slow manifold $\mathcal{M} = S_{\eps}$ and its accompanying linear fast fiber bundle $\mathcal{N}=\mathcal{N}_{\eps}$, as defined by Fenichel's theorem (Sec. \ref{sec:fenichel}). The characterizations \eqref{eq:manifold}--\eqref{eq:fiber} then state that $S_{\eps}$ is defined by the locus of points $B_{s\bot}(z)H(z) = 0$ where the components of the vector field $H(z)$ lying in the direction of the linear fast fibers vanish, and that the linear fast fibers $\mathcal{N}_p$ at basepoints $p \in S_{\eps}$ are spanned by the columns of the block component $A_f(p)$.

%

%
%
%The reader is referred to the appendix (and references therein) for a condensed introduction to the framework underlying the CSP method. Here, we focus on the details of the update step which are relevant for the proceeding section.\\

\subsection{CSP objects} The CSP iteration acts on a suitably initialised (point-dependent) basis matrix $A^{(0)}(z)$ and dual $B^{(0)}(z)$, producing a sequence $\{(A^{(j)}(z),B^{(j)}(z))\}_{j=0}^{\infty}$ of successively refined bases. We initialise with the pair 
\begin{eqnarray}
A^{(0)} &=& \begin{pmatrix}
A_f^{(0)} & A_s^{(0)}
\end{pmatrix} \label{eq:initialcond}\\
B^{(0)} &=& \begin{pmatrix}
B_{s\bot}^{(0)} \\ B_{f\bot}^{(0)}
\end{pmatrix}.\nn
\end{eqnarray}
The sequence $\{(A^{(j)}(z),B^{(j)}(z))\}_{j=0}^{\infty}$ in turn defines a sequence of updates $\{\Lambda^{(j)}\}_{j=0}^{\infty}$ to the $\Lambda$ operator:
\begin{eqnarray}
\Lambda^{(j)} &=& B^{(j)}[A^{(j)},H]. \label{eq:lambdaiter}
\end{eqnarray}
These approximate operators may be written in block components in analogy to \eqref{eq:lambda}. 
Motivated by the characterizations of the invariant manifold and fast fibers provided in Eqs. \eqref{eq:manifold}--\eqref{eq:fiber}, we define the following approximating objects: \\

\begin{definition}
~The {\rm CSP manifold of order $0$} is the level set
\begin{eqnarray}\label{eq:csp-k0}
\mathcal{K}^{(0)}  &=& \{(x,y) \in \mathbb{R}^n: B^{(0)}_{s\bot}(x,y,\eps) H(x,y,\eps) =  \mathbb{O}_{n-k,1}\}. \label{eq:cspmanifold0}
\end{eqnarray}
For integers $j \geq 1$, the {\rm CSP manifold of order} $j$, denoted $\mathcal{K}^{(j)}$, is
\begin{eqnarray}
\mathcal{K}^{(j)} &=& \{(x,y) \in \mathbb{R}^n: B^{(j)}_{s\bot}(x,\psi^{(j-1)}(x,\eps),\eps) H(x,y,\eps) = \mathbb{O}_{n-k,1}\}, \label{eq:cspmanifold} 
\end{eqnarray}
where $ y = \psi^{(j-1)}(x,\eps)$ is a graph of $\mathcal{K}^{(j-1)}$.  ~\\
\end{definition}

\begin{definition}
~For $j \geq 0$ and $p \in \mathcal{K}^{(j)}$, the {\rm CSP fiber of order} $j$ is the subspace
\begin{eqnarray}
\mathcal{L}^{(j)}(p)&=& \text{Col } A^{(j)}_f (p,\eps). \label{eq:cspfibers}
\end{eqnarray}
The {\rm CSP fiber bundle of order $j$} is the corresponding vector bundle
\begin{eqnarray}
\mathcal{L}^{(j)}&=& \bigcup_{p \in \mathcal{K}^{(j)}} \mathcal{L}^{(j)}_{\eps}(p). \label{eq:cspfiberbundle}
\end{eqnarray}
\end{definition}
\begin{remark}
~The convergence of the CSP manifolds and fiber bundles to the invariant slow manifold and fast fiber bundle of system \eqref{eq:master} depends on the choice of iteration step, as shown in Sec. \ref{sec:cspconvergence}. \\
\end{remark}

\subsection{One-step and two-step CSP updates} There are two commonly-used variants of the CSP iteration. We will only apply the two-step method in this paper, but it is instructive to introduce the simpler one-step method to clarify the relationship between these iterations and the CSP objects \eqref{eq:cspmanifold}--\eqref{eq:cspfibers}. 

Both methods use near-identity transformations in the update step, i.e. multiplication by matrices of the form $I\pm U$ and $I \pm L$, where $I$ is the identity matrix, and $U$ and $L$ are nilpotent matrices of the form

 \begin{eqnarray}
U^{(j)} &=& \begin{pmatrix}
\mathbb{O}_{n-k,n-k} & \tilde{U}^{(j)}\\
\mathbb{O}_{k,n-k} & \mathbb{O}_{k,k}
\end{pmatrix} \label{eq:lumatrices} \\
L^{(j)} &=& \begin{pmatrix}
\mathbb{O}_{n-k,n-k} & \mathbb{O}_{n-k,k} \\
\tilde{L}^{(j)} & \mathbb{O}_{k,k}
\end{pmatrix}. \nn
\end{eqnarray}
The block components $\tilde{U}^{(j)}$ and $\tilde{L}^{(j)}$,  of respective sizes $(n-k) \times k$  and $k \times (n-k)$, are defined by the constraint that the CPS manifolds and fibers converge asymptotically to the slow manifold and linear fast fibers; see Sec. \ref{sec:cspconvergence}. Near-identity update matrices are computationally easy to invert. In particular, we obtain efficient update rules for the dual basis $B^{(j)}$, as shown in the following two methods.

\subsubsection{One-step CSP method} \label{sec:csp1step}
The one-step CSP method is given by the iteration rule

\begin{eqnarray}
A^{(j+1)} &=& A^{(j)} (I- U^{(j)}) \label{eq:csp1}\\
&=& \begin{pmatrix} A_f^{(j)} & A_s^{(j)} - A_f^{(j)}\tilde{U}^{(j)} \end{pmatrix} \nn \\
B^{(j+1)} &=& (I+ U^{(j)}) B^{(j)} \nn\\
&=&  \begin{pmatrix} B_{s\bot}^{(j)} + \tilde{L}^{(j)} B_{f\bot}^{(j)} \\  B_{f\bot}^{(j)}\end{pmatrix}.\nn
\end{eqnarray}
The one-step CSP method only updates `half' of the basis. This update is sufficient if we are interested in computing only the CSP manifolds \eqref{eq:cspmanifold}. If we wish to approximate the fast fibers to the invariant manifold in tandem, we must update the remaining blocks of the basis and dual basis matrices. This update is provided by the two-step CSP method.\\

\subsubsection{Two-step CSP method} \label{sec:csp2step}
The two-step CSP method is given by the iteration rule
\begin{eqnarray}
A^{(j+1)} &=& A^{(j)} (I- U^{(j)})(I+L^{(j)}) \label{eq:csp2}\\
&=& \begin{pmatrix} A_f^{(j)}(I-\tilde{U}^{(j)}\tilde{L}^{(j)}) + A_s^{(j)} \tilde{L}^{(j)} & A_s^{(j)} - A_f^{(j)}\tilde{U}^{(j)}\end{pmatrix}\nn\\
B^{(j+1)} &=& (I-L^{(j)})(I+ U^{(j)}) B^{(j)}, \nn \\
&=& \begin{pmatrix}  B_{s\bot}^{(j)} + \tilde{U}^{(j)} B_{f\bot}^{(j)} \\ (I- \tilde{L}^{(j)} \tilde{U}^{(j)} ) B_{f\bot}^{(j)} - \tilde{L}^{(j)}B_{s\bot}^{(j)}\end{pmatrix}.\nn
\end{eqnarray}

%We define the {\it CSP manifold of order $j$}, $\mathcal{K}^{(j)}$, by the level set 
%
%\begin{eqnarray}
%B^{(j)}_1(x,\psi^{(j-1)}(x,\eps),\eps) H(x,y,\eps) &=& \mathbb{O}_{n-k,1}, \label{eq:cspmanifold}
%\end{eqnarray}

%We now turn to the {\it two-step} CSP method, which updates the remaining block column (resp. row) of the $A^{(j)}$ ($B^{(j)}$) basis matrices. The update step is
%
%\begin{eqnarray}
%A^{(j+1)} &=& A^{(j)} (I- U^{(j)})(I+L^{(j)})\\
%B^{(j+1)} &=& (I-L^{(j)})(I+ U^{(j)}) B^{(j)},
%\end{eqnarray}
%
%where we introduce the new near-identity transformation
%
%\begin{eqnarray}
%L^{(j)} &=& \begin{pmatrix}
%\mathbb{O}_{n-k,n-k} & \mathbb{O}_{n-k,k} \\
%L_{12}^{(j)} & \mathbb{O}_{k,k}
%\end{pmatrix}.
%\end{eqnarray}
%
%The two-step CSP method updates the first block column of $A^{(j)}$ and the second block column of $B^{(j)}$. We define the {\it CSP fiber of order $j$ at $z$}, $\mathcal{L}^{(j)}_{\eps}(z)$, by
%
%\begin{eqnarray}
%\mathcal{L}^{(j)}_{\eps}(z)&=& \text{span } \text{Col } A^{(j)}_1 (z,\eps).
%\end{eqnarray}
%
%
%The lower block is defined by the formula
%
%\begin{eqnarray}
%L_{12}^{(j)} &=& \Lambda_{21}^{(j)} (\Lambda_{11}^{(j)})^{-1}.
%\end{eqnarray}

\subsection{Convergence} \label{sec:cspconvergence}

Suppose that after $j$ iterates of the CSP two-step method, the CSP manifold of order $j$ defined by \eqref{eq:cspmanifold} is locally expressed as a graph $y = \psi^{(j)}(x,\eps)$ and then expanded as an asymptotic series in $\eps$:
\begin{eqnarray}
y &=& \psi^{(j)}(x,\eps) = \psi^{(j)}_0(x) + \eps \psi^{(j)}_1(x) +   \eps^2 \psi^{(j)}_2(x) + \cdots. \label{eq:cspmanifoldexp}
\end{eqnarray}
When the update terms in the near-identity transformations \eqref{eq:csp2} are defined appropriately, Kaper, Kaper, and Zagaris \cite{kkz2004a} demonstrated convergence of the CSP manifold to the invariant slow manifold of Fenichel's theory (as described in Sec. \ref{sec:fenichel}), in the following sense.~\\

\begin{theorem}[Convergence of CSP manifolds \cite{kkz2004a}] \label{thm:manifold}~Suppose the leading-order term of the graph of $\mathcal{K}^{(0)}$, denoted $y = \psi^{(0)}(x,\eps)$, agrees with the graph of the critical manifold $y = h_0(x)$ under a local choice of coordinates. Then using the one-step or two-step CSP update rule and for $j$ fixed, we have 
\begin{eqnarray*}
y = \psi^{(j)}(x,\eps) &=& \sum_{i=0}^j \eps^i h_i(x) + \mathcal{O}(\eps^{j+1})
\end{eqnarray*}

when $\eps > 0$ is sufficiently small, where $\psi^{(j)}$ is the graph of $\mathcal{K}^{(j)}$ with asymptotic series expansion \eqref{eq:cspmanifoldexp} and $y = h(x) = \sum_{j=0}^{\infty}h_j(x)$ is the asymptotic series expansion of the graph of $S_{\eps}$. 
\end{theorem}

This theorem states if the CSP two-step method is initialised appropriately, then for pairs of indices $0 \leq i \leq j$ we have $\psi^{(j)}_i(x) = h_i(x)$. This justifies the nomenclature `CSP manifold {\it of order $j$}': the set $\mathcal{K}^{(j)}$ agrees with $S_{\eps}$ up to order $j$ terms in the asymptotic expansion. 

An analogous convergence statement can be given for the CSP fiber bundle \eqref{eq:cspfiberbundle}: \\

\begin{theorem}[Convergence of CSP fibers \cite{kkz2004b}]\label{thm:fibers}~Fix $j\geq 0$ and let  $\mathcal{K}^{(j)}$ be given locally by the graph $y = \psi^{(j)}(x)$ for $x \in \mathbb{R}^k$ and $y \in \mathbb{R}^{n-k}$. Then given the basepoint $(x,\psi^{(j)}(x))$ on $\mathcal{K}^{(j)}$ (resp. $(x,h_{\eps}(x))$ on $S_{\eps}$), the asymptotic expansions of $\mathcal{L}_{\eps}^{(j)}(x,\psi^{(j)}(x))$ and $\mathcal{N}_{\eps}(x,h(x))$  agree up to and including terms of $\mathcal{O}(\eps^j)$. Thus, the CSP fast fiber bundle $\mathcal{L}^{(j)}$ is an $\mathcal{O}(\eps^j)$ approximation to the fast fiber bundle $\mathcal{N}_{\eps}$. 
\end{theorem}

Theorems \ref{thm:manifold}--\ref{thm:fibers} hold when the near-identity updates in \eqref{eq:lumatrices} are defined as follows:
\begin{eqnarray}
\tilde{U}^{(j)} &=& (\Lambda_{ff}^{(j)})^{-1} \Lambda_{fs}^{(j)} \label{eq:auxblock}\\
\tilde{L}^{(j)} &=& \Lambda_{sf}^{(j)} (\Lambda_{ff}^{(j)})^{-1}. \nn
\end{eqnarray}

With this choice of update step, the CSP two-step method simultaneously block-diagonalizes the CSP operator $\Lambda^{(j)}$ (as defined in \eqref{eq:lambdaiter}) in discrete steps as follows  \cite{kkz2004b}:\\

\begin{lemma} \label{lemma:lambdaconv}
For $j = 0,1,\cdots$, the CSP two-step method \eqref{eq:csp2} provides the asymptotic estimates
\begin{eqnarray*}
\Lambda^{(j)} &=& \begin{pmatrix}
\Lambda_{ff,0}+\mathcal{O}(\eps) & \mathcal{O}(\eps^j) \\
\mathcal{O}(\eps^j) &  \mathcal{O}(\eps)
\end{pmatrix},
\end{eqnarray*}

where $\Lambda^{(j)}$ is evaluated on $\mathcal{K}^{(j)}$. 
\end{lemma}

\section{Nonstandard CSP updates}\label{sec:csp3}

We now describe the CSP update step for nonstandard slow-fast systems \eqref{eq:master} satisfying Assumptions \ref{ass1}--\ref{assu1b}. The key is that bases for the fast and slow subspaces can be `read off' using the factorization $H(z,0)=N(z)f(z)$, providing a natural initial condition for the iteration. By Lemma \ref{lemma2}, the columns of $N(z)$ span the linear fast fiber with basepoint $z \in S$, and the columns of $Df(z)^{\top}$ span $T_z S^{\bot}$. Thus,
\begin{eqnarray}
A^{(0)} &=& \begin{pmatrix}
N & (Df^{\top})^{\bot}
\end{pmatrix}, \label{eq:initialize}
\end{eqnarray}
where the notation $P^{\bot}$ refers to a matrix whose columns form a basis to the subspace orthogonal to the column space of the matrix $P$. We define
\begin{eqnarray}
B^{(0)} &=& \begin{pmatrix}
(Df N)^{-1}Df \\
((N^{\bot})^{\top}(Df^{\top})^{\bot})^{-1} (N^{\bot})^{\top}
\end{pmatrix},\nonumber
\end{eqnarray}
where the matrix prefactors in both block rows (which are well-defined by Lemma \ref{lemma2}) normalize the block diagonal components in the product $B^{(0)}A^{(0)}$ to the identity, as required. 
With this initialization, the leading-order approximation $\mathcal{K}^{(0)}$ of the CSP manifold defined in \eqref{eq:csp-k0} can be computed:
\begin{eqnarray}
B^{(0)}_{s\bot} H &=& 0\\
(DfN)^{-1} Df (N f + \eps G) &=& 0 \nonumber\\
f + \eps (DfN)^{-1} (DfG) &=& 0 \nonumber\\
f &=& - \eps (DfN)^{-1} (DfG). \label{eq:k0}
\end{eqnarray}
This level set implicitly defines the graph $y = \psi^{(0)}(x,\eps)$. On this graph, we expand both sides of this equation in powers of $\eps$. The leading-order $O(\eps^0)$ coefficient is
\begin{eqnarray}\label{eq:leading-f}
f(x,\psi_0(x),0) &=& 0\,,
\end{eqnarray}
which matches the definition of the critical manifold $S$ as the leading-order part of the asymptotic series of $\mathcal{S}_{\eps}$, $f(x,h_0(x),0) = 0$; thus, $\psi^{(0)}_0(x) = h_0(x)$. 

The four block components of the operator $\Lambda^{(0)} = B^{(0)} [A^{(0)},H]$ can be computed explicitly in terms of the components $N,f,$ and $G$ in the vector field $H$; compare with \eqref{eq:lambda}:
\begin{eqnarray}
\Lambda_{ff}^{(0)} &=& B_{s\bot}^{(0)} [A_f^{(0)},H] \label{eq:lambdaff}\\
&=& (DfN)^{-1} Df(DH N - DN H)\nonumber\\ 
&=& (DfN)^{-1} Df ((NDf + \eps DG) N - DN (Nf + \eps G))\nonumber\\
&=& DfN - (DfN)^{-1} ( DfDN Nf)+ \eps (DfN)^{-1}Df [N,G],\nonumber
\end{eqnarray}
\begin{eqnarray}
\Lambda_{fs}^{(0)} &=& B_{s\bot}^{(0)} [A_s^{(0)},H] \label{eq:lambdafs}\\
&=& (DfN)^{-1} Df ((NDf + \eps DG) (Df^{\top})^{\bot}  - D((Df^{\top})^{\bot}) (Nf + \eps G))\nonumber\\
&=& - D((Df^{\top})^{\bot}) Nf  + \eps (DfN)^{-1} Df [(Df^{\top})^{\bot},G],\nonumber
\end{eqnarray}
\begin{eqnarray}
\Lambda_{sf}^{(0)} &=& B_{f\bot}^{(0)} [A_f^{(0)},H] \label{eq:lambdasf}\\
&=& ((N^{\bot})^{\top}(Df^{\top})^{\bot})^{-1} (N^{\bot})^{\top} (DH N - DN H)\nonumber\\
&=&  ((N^{\bot})^{\top}(Df^{\top})^{\bot})^{-1} (N^{\bot})^{\top} ((NDf + \eps DG) N - DN (Nf + \eps G))\nonumber\\
&=&  ((N^{\bot})^{\top}(Df^{\top})^{\bot})^{-1} (N^{\bot})^{\top} (-DN Nf + \eps [N,G])\nonumber,
\end{eqnarray}
%
%and
%
\begin{eqnarray}
\Lambda_{ss}^{(0)} &=& B_{f\bot}^{(0)} [A_s^{(0)},H] \label{eq:lambdass}\\
&=& ((N^{\bot})^{\top}(Df^{\top})^{\bot})^{-1} (N^{\bot})^{\top} ((NDf + \eps DG) (Df^{\top})^{\bot} \nn\\
&&  - D((Df^{\top})^{\bot}) (Nf + \eps G))\nonumber\\
&=&  ((N^{\bot})^{\top}(Df^{\top})^{\bot})^{-1} (N^{\bot})^{\top}(-D((Df^{\top})^{\bot}) Nf +\eps[ (Df^{\top})^{\bot},G])\nonumber
\end{eqnarray}
The vector function $f$ is $\mathcal{O}(\eps)$ on $\mathcal{K}^{(0)}$ by \eqref{eq:k0}; therefore, the latter three block components are of at most $\mathcal{O}(\eps)$ as well. This in turn implies that the update blocks $\tilde{U}^{(0)}$ and $\tilde{L}^{(0)}$ defined in \eqref{eq:auxblock} are both $\mathcal{O}(\eps)$ on this set. 
After one application of the CSP two-step update \eqref{eq:csp2}, the updated bases $A^{(1)}$ and $B^{(1)}$ are
\begin{eqnarray}
A^{(1)} &=& \begin{pmatrix}
A_f^{(0)} + A_s^{(0)} \tilde{L}^{(0)}-A_f^{(0)}\tilde{U}^{(0)}\tilde{L}^{(0)} &
-A_f^{(0)} \tilde{U}^{(0)} + A_s^{(0)}
\end{pmatrix}\\
B^{(1)} &=& \begin{pmatrix}
 B_{s\bot}^{(0)} + \tilde{U}^{(0)} B_{f\bot}^{(0)} \\
 B_{f\bot}^{(0)} - \tilde{L}^{(0)} B_{s\bot}^{(0)} - \tilde{L}^{(0)}\tilde{U}^{(0)} B_{f\bot}^{(0)}
\end{pmatrix}.
\end{eqnarray}
Without writing out the update in detail, it suffices to note that the block components $\tilde{L}^{(0)}$ and $\tilde{U}^{(0)}$ each introduce $\mathcal{O}(\eps^1)$ perturbations of the initial block columns of $A^{(0)}$ (respectively block rows of $B^{(0)}$). The first fast fiber update is explored in greater detail using projection maps, in Sec. \ref{sec:proj}. Subsequent CSP iterations will result in higher-order perturbations. 

The first update of the CSP manifold $\mathcal{K}^{(1)}$ is defined as follows; see \eqref{eq:cspmanifold}:
\begin{eqnarray}
B_{s\bot}^{(1)}(x,\psi^{(0)}(x,\eps),\eps) H(x,y,\eps) &=& 0\\
(D\!f_0 N_0)^{-1} D\!f_0 (Nf+\eps G) + \Lambda_{fs,0} B_{f\bot,0} (Nf + \eps G) &=& 0.\nonumber
\end{eqnarray}
In this equation and the following, we use the additional subscript '$0$' to refer to those quantities that are computed on points $(x,\psi^{(0)}(x))\in \mathcal{K}^{(0)}$. We compute the order $O(\eps^1)$ coefficient of the asymptotic series assuming the graph form $y = \psi^{(1)}(x,\eps)$ where the leading order coefficient is defined by $f_0 = 0\,$; see \eqref{eq:leading-f}:
\begin{eqnarray}
(Df_0 N_0)^{-1} Df_0 (N_0f_0 + \eps (N_0 D_y f_0 \psi^{(1)}_1 + G_0)) + &&\\
 \Lambda_{fs,0} B_{f\bot,0} (N_0f_0 + \eps ( N_0 D_y f_0  \psi^{(1)}_1 + G_0) + \mathcal{O}(\eps^2)) &=& 0\nonumber\\
 \eps D_y f_0  \psi^{(1)}_1 + \eps (Df_0 N_0)^{-1} Df_0 G_0 + O(\eps^2) &=& 0\,,\nonumber
\end{eqnarray}
which gives
\begin{eqnarray}
 \psi^{(1)}_1 &=& -(D_y f_0)^{-1}(Df_0 N_0)^{-1} (Df_0 G_0). \label{eq:firstorder}
\end{eqnarray}
This result matches a recent result due to Wechselberger \cite{wechselberger2019} for the first-order correction of $\mathcal{S}_{\eps}$ in the nonstandard case.\\

%\newpage
\section{Examples}\label{sec:eg}

We now demonstrate the two-step CSP method in four nonstandard examples. Our list of systems increases in complexity: first we consider the case where the slow manifold is equal to the critical manifold and the fast fiber bundle remains unchanged. Second, we consider the case where the slow manifold is equal to the critical manifold but the fast fiber bundle perturbs from the $\eps = 0$ case. Third, we consider a system where both the manifold and fiber updates are nontrivial. Finally, we consider the four-dimensional system \eqref{eq:slime2} discussed in the introduction. \\

\subsection{Trivial updates of the slow manifold and fast fibers}\label{sec:circle}

Consider the following two-dimensional nonstandard system:
 
 \begin{equation} \label{eq:circ}
\begin{pmatrix}
x' \\ y'
\end{pmatrix} = 
\begin{pmatrix}
x \\ y
\end{pmatrix} 
(1-x^2-y^2) + \eps
\begin{pmatrix}
-y \\
x
\end{pmatrix}.
\end{equation}
%
\begin{comment}
 \begin{eqnarray}
x' &=& x(1-(x^2+y^2)) - \eps y \label{eq:circ}\\
y' &=& y(1-(x^2+y^2)) + \eps x.\nonumber
\end{eqnarray}

Defining $z = (x,y)$ and 
\begin{eqnarray}
N(x,y) &=& \begin{pmatrix}
x \\y
\end{pmatrix}\\
f(x,y) &=& 1-(x^2+y^2)\nonumber\\
G(x,y,\eps) &=& \begin{pmatrix}
-y \\ x
\end{pmatrix},\nonumber
\end{eqnarray}

the system can be written in the factorized form \eqref{eq:master}.
\end{comment}
%
The critical manifold $S = \{x^2 + y^2 = 1\}$ is a circle. The Jacobian evaluated along $S$ is
\begin{equation}
DH|_{S} =  \left. \begin{pmatrix}
-2x^2 & -2xy\\ -2xy & -2y^2 
\end{pmatrix}
\right|_S\,,
\end{equation}
with one trivial zero eigenvalue and one negative eigenvalue $\lambda = -2$ for all $(x,y)\in S$, implying that the critical manifold $S$ is attracting and normally hyperbolic. Thus, Assumptions \ref{ass1}--\ref{assu1b} are satisfied. Note that the origin is an isolated singularity of $N(z)$, where $f(z) \neq 0$, but we do not consider this singularity in the CSP iteration (see Remark \ref{remark:isolated}).

\begin{figure}[t]
  \centering
    \includegraphics[width=0.95\textwidth]{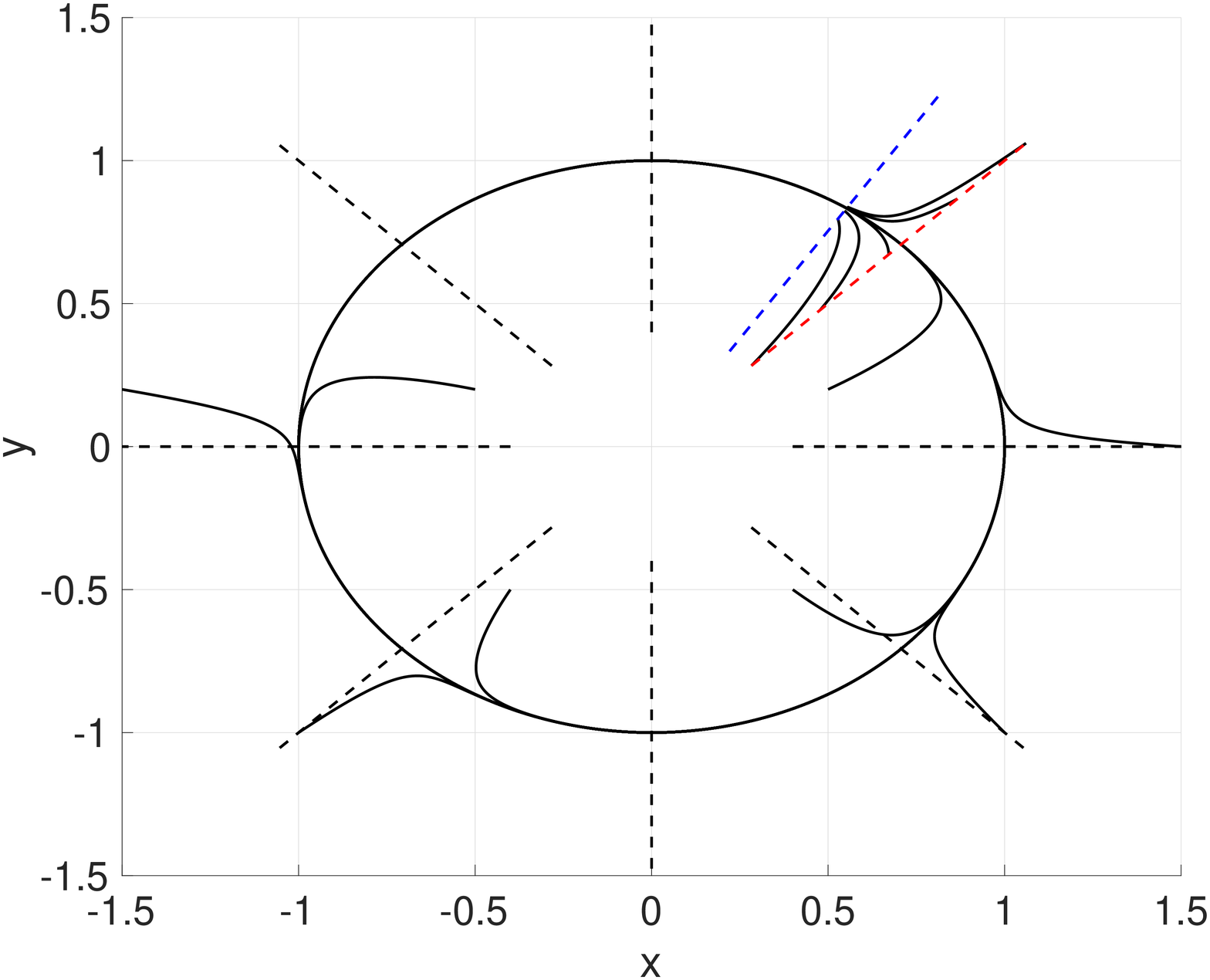}
      \caption{Several trajectories (black solid curves) plotted for system \eqref{eq:circ} (with $\eps = 0.1$). Portions of the fast fibers given in \eqref{eq:circlebundle} are given by dashed lines. The forward invariance of the family is also illustrated: trajectories on the red fiber $\mathcal{F}_{\eps}(p(\pi/4))$ are flowed forward for a time $t = 2$. These trajectories all end on the blue fiber $\mathcal{F}_{\eps}(p(\pi/4+2\eps))$. Trajectories computed using a Dormand-Prince ODE solver in MATLAB R2018a \cite{matlab}.}
      \label{fig:circle}
\end{figure}

When $\eps$ is small and positive, there is a clear separation between slow motion along $S$ versus fast motion toward $S$. The dynamics is made obvious if we write the system in polar coordinates and plot a graph of representative trajectories (Fig. \ref{fig:circle}):
\begin{eqnarray}
r' &=& r(1-r^2) \nn\\
\theta' &=& \eps. \label{eq:polar}
\end{eqnarray}

The coordinate representation \eqref{eq:polar} decomposes the system into two independent ODEs. The invariant manifold $r=1$, which is independent of $\eps$, can be read off from the first equation of \eqref{eq:polar}. Thus, $S_{\eps} = S$. We can also determine that the (nonlinear) fast fiber bundle $\mathcal{F}_{\eps}$ is foliated by the family of rays extending from the origin, with basepoints on the circle: 
\begin{eqnarray}
\mathcal{F}_{\eps}(p(\theta_0)) &=& \{(r,\theta): 0<r<\infty,~\theta=\theta_0\}, \hspace{0.5cm} p(\theta_0) = (1,\theta_0) \in S_{\eps} \nn \\
\mathcal{F}_{\eps} &=& \bigcup_{\theta_0 \in [0,2\pi)} \mathcal{F}_{\eps}(p(\theta_0)). \label{eq:circlebundle}
\end{eqnarray}

Note that the nonlinear and linear fast fiber bundles coincide in this example: $\mathcal{F}_{\eps} = \mathcal{N}_{\eps}$. The fast fiber bundle is invariant as a family: for $t\in \mathbb{R}$, the fiber $\mathcal{F}_{\eps}(p(\theta_0))$ is mapped to $\mathcal{F}_{\eps}(p(\theta_0+\eps t))$ under the time $t$ flow map. Furthermore, exponential contraction of the rays toward $S_{\eps}$ is governed by the equation $r' = r(1-r^2)$. Compare these facts about the nonlinear fast fibers with Fenichel's theorem, Sec. \ref{sec:fenichel}. This behavior is similarly independent of the value of $\eps > 0$.

We initialize the CSP two-step method using \eqref{eq:initialize}: 
\begin{eqnarray}
A^{(0)} &=& \begin{pmatrix}
N & (Df^{\top})^{\bot}
\end{pmatrix} \nn\\
&=& \begin{pmatrix}
x  & 2y\\ 
y & -2x
\end{pmatrix} \label{eq:circlebasis1}\\
B^{(0)} &=& \frac{1}{2(x^2+y^2)}\begin{pmatrix}
2x & 2y \\
y & -x
\end{pmatrix}. \label{eq:circlebasis2}
\end{eqnarray}
We begin by computing the initial CSP manifold $\mathcal{K}^{(0)}$  (see \eqref{eq:cspmanifold0}). We compute
\begin{eqnarray}
B^{(0)}_{s\bot}H &=& \frac{1}{2(x^2+y^2)}\begin{pmatrix}
2x & 2y
\end{pmatrix} \left( \begin{pmatrix}
x \\ y
\end{pmatrix} (1-x^2-y^2) + \eps \begin{pmatrix}
-y \\ x
\end{pmatrix}\right) \\
&=& 1-x^2-y^2 = 0, \nonumber
\end{eqnarray}
and so 
\begin{eqnarray*}
\mathcal{K}^{(0)} =  \{(x,y):x^2 + y^2 = 1\} = S.
\end{eqnarray*}
We have
\begin{eqnarray*}
DA_f|_S &=& \begin{pmatrix}
1 & 0 \\ 0 & 1
\end{pmatrix} \\
DA_s|_S &=& \begin{pmatrix}
0 & 2 \\ -2 & 0\end{pmatrix}\\
DH|_S &=& \begin{pmatrix}
-2x^2 & -2xy - \eps \\
-2xy + \eps & -2y^2
\end{pmatrix}.
\end{eqnarray*}
With this information, the initial $\Lambda^{(0)}$ matrix \eqref{eq:lambdaiter} becomes
\begin{eqnarray}
\Lambda^{(0)}|_S &=& \begin{pmatrix}
-2 &  0 \\ 0 & 0
\end{pmatrix}. \label{eq:circlelambda}
\end{eqnarray}
Using \eqref{eq:auxblock}, we find that  $\tilde{U}^{(0)} = 0$ and $\tilde{L}^{(0)} = 0$.

The triviality of the updates implies that the CSP manifolds and fibers  (\eqref{eq:cspmanifold}--\eqref{eq:cspfibers})  are
\begin{eqnarray*}
\mathcal{K}^{(j)} &=& \{(x,y) \in \mathbb{R}^2: x^2 + y^2 = 1\} = S\\
\mathcal{L}^{(j)}(p) &=& \{cp: c\in \mathbb{R}\}, \hspace{0.5cm}p \in S\\
\end{eqnarray*}
for $j = 0,1,2,\cdots$.

\subsection{Trivial updates of the slow manifold; nontrivial updates of the fast fibers}\label{sec:modcircle}

We now consider a variant of the previous system, where we modify the first component of $G(x,y)$: 
\begin{equation} \label{eq:circ2}
\begin{pmatrix}
x' \\ y'
\end{pmatrix} = \begin{pmatrix}
x \\ y
\end{pmatrix} (1-x^2-y^2) + \eps
\begin{pmatrix}
-y + y^2 (1-x^2-y^2)\\
x
\end{pmatrix}.
\end{equation}
We still have $S = \{(x,y): x^2 + y^2 = 1\}$, and it is easy to check that this is also the invariant slow manifold for $\eps > 0$: for $p=(x,y) \in S$ we have
\begin{eqnarray}
H(x,y,\eps) &=&   N(x,y,\eps)  f(x,y) + \eps G(x,y,\eps)\\
 &=& \eps G(x,y,\eps) \nonumber\\ 
&=& \begin{pmatrix}
- \eps y \\ \eps x
\end{pmatrix},
\end{eqnarray}
so that $H(p,\eps)\in T_p S$.  Thus, $S = S_{\eps}$ as before. In this modified system, however, the linear fast fiber $\mathcal{N}_{\eps,p}$ will not be orthogonal to $T_p S_{\eps}$ at points $p \in S_{\eps}$. 

The CSP step is initialized as usual with \eqref{eq:initialize}. The basis matrices are identical to \eqref{eq:circlebasis1}--\eqref{eq:circlebasis2} since we have not modified $N(z)$ or $f(z)$. From \eqref{eq:lambdaiter}, the initial operator $\Lambda^{(0)}$ is (compare \eqref{eq:circlelambda})
\begin{eqnarray}
\Lambda^{(0)}|_{S} &=& \begin{pmatrix}
-2 - 2\eps xy^2 & 0 \\ 
\eps y^3 & 0 
\end{pmatrix}.
\end{eqnarray}
By definition \eqref{eq:auxblock}, the update term $\tilde{U}^{(0)}$ is trivial since $\Lambda_{fs}^{(0)} = 0$; on the other hand, the update quantity $\tilde{L}^{(0)}$ may give a nontrivial first-order correction to the fast fiber. The first column (resp. first row) of the updated basis matrix $A^{(1)}$ (resp. $B^{(1)}$) are
\begin{eqnarray}
A^{(1)}_f|_S  &=& \begin{pmatrix}
x + y^4 \eps \\
y - xy^3 \eps 
\end{pmatrix} + \mathcal{O}(\eps^2) \nonumber\\
B^{(1)}_{s\bot}|_S &=& \begin{pmatrix}
x & y
\end{pmatrix}. 
\end{eqnarray}

Applying the definitions \eqref{eq:cspmanifold}--\eqref{eq:cspfibers}, the updated CSP objects are 
\begin{eqnarray*}
\mathcal{K}^{(1)} &=& \{B^{(1)}_{s\bot}H = 0 \}\\
&=& \{(x,y) \in \mathbb{R}^2: x^2 + y^2 = 1\}\\
\mathcal{L}^{(1)}(x,y(x),\eps) &=& \left\{ c \left(\begin{pmatrix}
x + y(x)^4 \eps \\
y - xy(x)^3 \eps 
\end{pmatrix} + \mathcal{O}(\eps^2) \right): c \in \mathbb{R}  \right\}, \hspace{0.5cm}, (x,y(x)) \in S. 
\end{eqnarray*}

Here, the function $y(x)$ refers to a graph of $x^2 + y^2 = 1$ containing the chosen basepoint $(x,y(x))$. Graphs of the form $x(y)$ can also be used. 

\subsection{Parabolic critical manifold}
Given the system
\begin{eqnarray}
\begin{pmatrix}
x' \\ y'
\end{pmatrix} &=& N(x,y) f(x,y) + \eps G(x,y,\eps)\\
&=& 
\begin{pmatrix}
-2x \\ -y
\end{pmatrix} (x^2 + y - 1)  + \eps \begin{pmatrix}
2 \\ -x + \eps
\end{pmatrix}.\nonumber
\end{eqnarray}
 The critical manifold of the system is defined as the level set $f = 0$, which in this case can be solved globally for the graph $y = 1 - x^2$. This critical manifold is globally attracting and normally hyperbolic---the Jacobian along $S$ has one trivial zero eigenvalue and another eigenvalue $-(3x^2+1) < 0$.

The CSP method is initialized with the following basis and dual as given in \eqref{eq:initialize}:
\begin{eqnarray}
A^{(0)} &=& \begin{pmatrix}
-2x & -1\\ -y & 2x
\end{pmatrix}\\
B^{(0)} &=& \frac{1}{4x^2 + y}\begin{pmatrix}
-2x & -1 \\
-y & 2x
\end{pmatrix}.
\end{eqnarray}

The initial CSP manifold \eqref{eq:cspmanifold0} written to $\mathcal{O}(\eps^{0})$ order is
\begin{eqnarray*}
\mathcal{K}^{(0)} = \{ (x,y): y = 1-x^2 + \mathcal{O}(\eps)\}.
\end{eqnarray*}

Define the auxiliary function $g(x) = 3x^2 + 1$. The operator $\Lambda^{(0)}$ (as defined in \eqref{eq:lambdaiter}) computed to linear order is
\begin{eqnarray*}
\Lambda^{(0)}|_{\mathcal{K}^{(0)}} &=& \begin{pmatrix}
-g(x) - \frac{12x}{g(x)}\eps & - 3\frac{3x^2-1}{g(x)^2} \eps\\
\frac{2(3x^2-2)}{g(x)}\eps & -\frac{12x}{g(x)^2}\eps,
\end{pmatrix} + \mathcal{O}(\eps^2).
\end{eqnarray*}
After one application of the two-step CSP method \eqref{eq:csp2}, we obtain the updated basis matrix
\begin{eqnarray*}
A^{(1)}|_{\mathcal{K}^{(0)}}&=& \begin{pmatrix}
A_f^{(1)} & A_s^{(1)}
\end{pmatrix}\\
&=& \begin{pmatrix}
-2x + \frac{2(3x^2 - 2)}{g(x)^2}  \eps + & -1 + \frac{6x(3x^2 - 1)}{g(x)^3}\eps \\
(x^2 - 1) - \frac{4(3x^3-2x)}{g(x)^2} \eps  & -2x - \frac{3(3x^4 - 4x^2 + 1)}{g(x)^3}\eps
\end{pmatrix} + \mathcal{O}(\eps^2)
\end{eqnarray*}

and dual $B^{(1)}|_{\mathcal{K}^{(0)}}$.

This gives us the first update of the CSP objects using the definitions  \eqref{eq:cspmanifold}--\eqref{eq:cspfibers}:
\begin{eqnarray*}
\mathcal{K}^{(1)} &=& \{(x,y) \in \mathbb{R}^2: y = \psi^{(1)}(x) = 1-x^2 + \frac{3x}{g(x)}\eps + \mathcal{O}(\eps^2) \}\\
\mathcal{L}^{(1)}(x,\psi^{(1)}(x)) &=& \text{span Col}(A^{(1)}_f(x,\psi^{1}(x)))\\
&=&  \left\{ c\left( \begin{pmatrix}
-2x + \frac{2(3x^2 - 2)}{g(x)^2} \eps \\
(x^2 - 1) - \frac{4(3x^3-2x)}{g(x)^2} \eps 
\end{pmatrix} + \mathcal{O}(\eps^2) \right): c \in \mathbb{R} \right\}.
\end{eqnarray*}

The $\Lambda^{(1)}$ update \eqref{eq:lambdaiter} is given by
\begin{eqnarray*}
\Lambda^{(1)}|_{\mathcal{K}^{(1)}} &=& \begin{pmatrix}
-g(x) - \frac{12x}{g(x)}\eps & 0 \\
0  & -\frac{12x}{g(x)^2}\eps,
\end{pmatrix} + \mathcal{O} (\eps^2).
\end{eqnarray*}
Observe that the off-diagonal elements are now 0 (modulo nonzero $\mathcal{O}(\eps^2)$ terms) while the diagonal terms are unchanged up to order $\eps$. As we continue to apply the CSP algorithm, the off-diagonal elements can be made to vanish modulo arbitrarily high orders; see Lemma \ref{lemma:lambdaconv}. The order we used in the initial basis matrix implies that the (1,1) term governs the fast dynamics and (2,2) term governs the slow dynamics in the decoupled system.

%A^{(1)}(x,\psi^{(1)}(x),\eps) 
%& -1 + \frac{6x(3x^2 - 1)}{g(x)^3}\eps
%& -2x - \frac{3(3x^4 - 4x^2 + 1)}{g(x)^3}\eps
%

\subsection{A four-dimensional system}

Recall system \eqref{eq:slime2}:
\begin{eqnarray}
\begin{pmatrix}
p' \\ d' \\ r' \\ b'
\end{pmatrix} &=& \begin{pmatrix}
2 & 2 \\
0 & 1 \\
1 & 0 \\
-1 & 0
\end{pmatrix} \begin{pmatrix}
-k_5 rp^2 + k_{-5}b\\
-k_4 dp^2 + k_{-4}(c-d-r-b)
\end{pmatrix} + \eps \begin{pmatrix}
k_3 - k_{-3}p + k_2 Sb \\
-k_1 d + k_{-1}r \\
k_1 d - k_{-1} r\\
0
\end{pmatrix}, \label{eq:slime3}
\end{eqnarray}
for parameters $k_i,S,c > 0$, together with the perturbation parameter $\eps>0$. This system has the factorized form $H(z) = N(z)f(z) + \eps G(z,\eps)$. Using \eqref{eq:initialize}, the initial bases are given by
\begin{eqnarray*}
A^{(0)} &=& \begin{pmatrix}
2 & 2 & \frac{k_{-5}}{2k_5 pr} & -\frac{p}{2r}\\
0 & 1 & -\frac{k_4 k_{-5}d + k_{-4}k_5 r}{k_5 (k_4 p^2 + k_{-4}) r} & -\frac{k_{-4}r - k_4 d p^2}{ (k_4 p^2 + k_{-4}) r}\\
1 & 0 & 0 & 1 \\
-1 & 0 & 1 & 0
\end{pmatrix}
\end{eqnarray*}
and the dual by $B^{(0)} = (A^{(0)})^{-1}$, whose first two rows are given by the $2\times 4$ matrix $(Df N)^{-1} Df$. 

As before, the leading-order term of the graph of $\mathcal{K}^{(0)}$ (defined in \eqref{eq:cspmanifold0}) will match the graph of $\{f = 0\}$:
\begin{eqnarray*}
\mathcal{K}^{(0)} &=& \left\{(p,d,r,b): (r,b) = \psi^{(0)}(p,d) =  \left( \frac{k_{-5}(-k_4 dp^2 + k_{-4} (c-d))}{k_4 ( k_5 p^2 + k_{-5})} ,\frac{k_5 p^2 (-k_4 dp^2 + k_{-4} (c-d))}{k_4 ( k_5 p^2 + k_{-5})}\right)  + \mathcal{O}(\eps) \right\}.
\end{eqnarray*}

Applying the CSP manifold definition \eqref{eq:cspmanifold}, we have
\begin{eqnarray*}
\mathcal{K}^{(1)} = \{(p,d,r,b):(r,b) &=& \psi^{(1)}_0(p,d) + \eps \psi^{(1)}_1(p,d) + \mathcal{O}(\eps^2)\},
\end{eqnarray*}
where $\psi^{(1)}_0(p,d) = \psi^{(0)}_0(p,d)$ and the first-order correction $ \psi^{(1)}_1(p,d)$ can be factored as
\begin{eqnarray*}
 \psi^{(1)}_1(p,d) &=& \frac{1}{\gamma(p)}(\beta_{2}(p) d^2 + \beta_{1}(p) d + \beta_0(p), \delta_{2}(p)d^2 + \delta_{1}(p)d+\delta_0(p)),
\end{eqnarray*}
for the seven polynomials $\beta_{i}, \delta_{i},\gamma$ in $p$, defined as follows:
\begin{eqnarray*}
\gamma(p) &=& (k_{-4}^2 (k_{-5} + k_{5} p^2)^2 (-4 d k_{4}^2 k_{5} k_{-5} p^5 +  k_{-4}^2 (k_{-5}^2 + k_{5}^2 p^4 + 2 k_{5} k_{-5} p (2 c - 2 d + p)) + \\
   && k_{4} k_{-4} p (4 d (k_{-5}^2 + k_{5}^2 p^4) +  p (k_{-5}^2 + k_{5}^2 p^4 + 2 k_{5} k_{-5} p (2 c + p)))) 
   \end{eqnarray*}   
   \begin{eqnarray*}
\beta_{2}(p) &=& k_2k_4k_{-1}k_1 p (-2 k_2 k_4 k_5^3 k_{-4} k_{-5} (k_4 p^2 + k_{-4}) S p^6 -  4 k_4 k_5^2 k_{-1} k_{-4} k_{-5}^2 (k_4 p^2 + k_{-4}) p^4 + \\
 &&  4 k_4 k_5^2 k_{-1} k_{-4}^2 k_{-5} (k_4 p^2 + k_{-4}) p^4 +  4 k_4^2 k_5 k_5 k_{-5}^2 (k_{-5} - k_{-4}) (k_4 p^2 + k_{-4}) p^4 + \\
  && 4 k_1 k_4 k_5 k_{-4} (k_{-4} - k_{-5}) (k_5 k_{-4} - k_4 k_{-5}) (k_5 p^2 + k_{-5}) p^4 - 2 k_2 k_4 k_5 k_{-4} k_{-5}^3 (k_4 p^2 + k_{-4}) S p^2 + \\
  && 2 k_3 k_5^2 k_{-4} k_{-5} (k_4 p^2 + k_{-5}) (k_4^2 p^4 + 2 k_{-4} (k_{-4} - k_{-5}) p^2 +    k_{-4}^2) S p^2) 
  \end{eqnarray*} 
  \begin{eqnarray*}
 \beta_{1}(p) &=& p (-2 k_{4} k_{5}^3 k_{-3} k_{-4}^2 k_{-5} p^7 + 2 c k_{2} k_{4} k_{5}^3 k_{-4}^2 k_{-5} S p^6 -   k_{4} k_{5}^2 k_{-1} k_{-4} k_{-5}^2 (k_{4} p^2 + k_{-4}) p^5 + \\
   &&    k_{4} k_{5}^2 k_{-1} k_{-4}^2 k_{-5} (k_{4} p^2 + k_{-4}) p^5 +     k_{1} k_{4} k_{5}^2 k_{-4}^2 (k_{-4} - k_{-5}) (k_{5} p^2 + k_{-5}) p^5 + \\
   &&    4 c k_{4} k_{5}^2 (k_{1} - k_{-1}) k_{-4}^2 (k_{-4} - k_{-5}) k_{-5} p^4 -     8 c k_{4}^2 k_{5} k_{-1} k_{-4} k_{-5}^2 (k_{-5} - k_{-4}) p^4 + \\
   &&    k_{4} k_{5} k_{-1} k_{-4} (k_{-4} - 2 k_{-5}) k_{-5}^2 (k_{4} p^2 + k_{-4}) p^3 +    k_{5}^2 k_{-1} k_{-4}^3 k_{-5} (k_{4} p^2 + k_{-4}) p^3 + \\
   &&    k_{1} k_{5} k_{-4}^2 (k_{5} k_{-4}^2 + k_{4} (k_{-4} - 2 k_{-5}) k_{-5}) (k_{5} p^2 + k_{-5}) p^3 +  2 k_{5}^2 k_{-3} k_{-4}^2 k_{-5} (k_{4}^2 p^4 + k_{4} (2 k_{-4} - 3 k_{-5}) p^2 + \\
   &&       k_{-4}^2) p^3 - 4 c k_{4} k_{5} k_{-1} k_{-4}^2 k_{-5}^2 (k_{-5} - k_{-4}) p^2 -     4 c k_{2} k_{5}^2 k_{-4}^2 k_{-5} (k_{4}^2 p^4 + k_{4} (2 k_{-4} - k_{-5}) p^2 + \\
   &&       k_{-4}^2) S p^2 +     2 c k_{4} k_{5} k_{-4}^2 k_{-5}^2 (2 k_{1} k_{-4} - 2 k_{1} k_{-5} + k_{2} k_{-5} S) p^2 + \\
   &&    k_{5} k_{-1} k_{-4}^3 k_{-5}^2 (k_{4} p^2 + k_{-4}) p -     k_{1} k_{4} k_{-4}^2 k_{-5}^3 (k_{5} p^2 + k_{-5}) p + \\
    &&   k_{1} k_{5} k_{-4}^4 k_{-5} (k_{5} p^2 + k_{-5}) p +    2 k_{5} k_{-3} k_{-4}^2 k_{-5}^2 (k_{4}^2 p^4 + k_{4} (2 k_{-4} - 3 k_{-5}) p^2 + \\
   &&       k_{-4}^2) p - k_{4} k_{-4} k_{-5}^4 (2 k_{-3} k_{-4} + k_{-1} (k_{4} p^2 + k_{-4})) p +    2 k3 k_{-4}^2 k_{-5} (k_{5} p^2 + k_{-5}) (k_{4} k_{5} (k_{5} - k_{4}) p^4 + \\
   &&       2 k_{4} k_{5} (k_{-5} - k_{-4}) p^2 - k_{5} k_{-4}^2 + k_{4} k_{-5}^2))
   \end{eqnarray*}
   \begin{eqnarray*}
    \beta_{0}(p) &=& p (c k_{4} k_{5}^2 k_{-4}^2 k_{-5} ((-k_{-1} - 2 k_{-3}) k_{-4} + k_{-1} k_{-5}) p^5 - 
   c k_{4} k_{5} k_{-1} k_{-4}^2 (k_{-4} - 2 k_{-5}) k_{-5}^2 p^3 + \\
  && c k_{5}^2 (-k_{-1} - 2 k_{-3}) k_{-4}^4 k_{-5} p^3 + 
   2 c^2 k_{5} k_{-4}^2 k_{-5} (2 k_{4} k_{-1} k_{-5} (k_{-5} - k_{-4}) + \\
  &&    k_{2} k_{5} k_{-4} (k_{4} p^2 + k_{-4}) S) p^2 + c k_{4} k_{-1} k_{-4}^2 k_{-5}^4 p - 
   c k_{5} k_{-1} k_{-4}^4 k_{-5}^2 p - 2 c k_{5} k_{-3} k_{-4}^3 k_{-5}^2 (k_{4} p^2 + k_{-4}) p + \\
 &&  2 c k_{3} k_{5} k_{-4}^3 k_{-5} (k_{4} p^2 + k_{-4}) (k_{5} p^2 + k_{-5}))
 \end{eqnarray*}
 \begin{eqnarray*}
 \delta_{2}(p) &=& -k_{5} p (-4 k_{4}^2 k_{5} k_{-1} k_{-5}^2 (k_{4} p^2 + k_{-4}) p^6 + 
   4 k_{4} k_{5}^2 k_{-1} k_{-4} k_{-5} (k_{4} p^2 + k_{-4}) p^6 +  4 k_{1} k_{4} k_{5}^2 k_{-4}^2 (k_{5} p^2 + k_{-5}) p^6\\
 &&  - 
   4 k_{1} k_{4}^2 k_{5} k_{-4} k_{-5} (k_{5} p^2 + k_{-5}) p^6 - 
   4 k_{4}^2 k_{-1} k_{-4} k_{-5}^2 (k_{4} p^2 + k_{-4}) p^4 + \\
  && 4 k_{4} k_{5} k_{-1} k_{-4}^2 k_{-5} (k_{4} p^2 + k_{-4}) p^4 + 
   4 k_{1} k_{4} k_{5} k_{-4}^3 (k_{5} p^2 + k_{-5}) p^4 - 
   4 k_{1} k_{4}^2 k_{-4}^2 k_{-5} (k_{5} p^2 + k_{-5}) p^4 + \\
   &&    2 k_{2} k_{5} k_{-4} (k_{4} p^2 + k_{-4}) (k_{4} k_{5}^2 p^6 + k_{4} (k_{4} + 2 k_{5}) k_{-5} p^4 + 
      k_{4} k_{-5} (2 k_{-4} + k_{-5}) p^2 + k_{-4}^2 k_{-5}) S p^2)
 \end{eqnarray*}
 \begin{eqnarray*}
 \delta_{1}(p) &=& -k_{5} p (k_{4} k_{5}^2 k_{-1} k_{-4} k_{-5} (k_{4} p^2 + k_{-4}) p^7 + 
   k_{1} k_{4} k_{5}^2 k_{-4}^2 (k_{5} p^2 + k_{-5}) p^7 + \\ &&
   8 c k_{4}^2 k_{5} k_{-1} k_{-4} k_{-5}^2 p^6 - 4 c k_{4} k_{5}^2 k_{-1} k_{-4}^2 k_{-5} p^6 + 
   2 k_{4} k_{5} k_{-1} k_{-4} k_{-5}^2 (k_{4} p^2 + k_{-4}) p^5 + \\ &&
   k_{4} k_{5} k_{-1} k_{-4}^2 k_{-5} (k_{4} p^2 + k_{-4}) p^5 + 
   k_{1} k_{4} k_{5} k_{-4}^3 (k_{5} p^2 + k_{-5}) p^5 + \\ &&
   2 k_{1} k_{4} k_{5} k_{-4}^2 k_{-5} (k_{5} p^2 + k_{-5}) p^5 + 
   4 c k_{4} (2 k_{4} + k_{5}) k_{-1} k_{-4}^2 k_{-5}^2 p^4 - \\ &&
   4 c k_{4} k_{5} k_{-1} k_{-4}^3 k_{-5} p^4 + k_{4} k_{-1} k_{-4} k_{-5}^3 (k_{4} p^2 + k_{-4}) p^3 +\\ &&
    k_{4} k_{-1} k_{-4}^2 k_{-5}^2 (k_{4} p^2 + k_{-4}) p^3 + 
   k_{1} k_{4} k_{-4}^2 k_{-5}^2 (k_{5} p^2 + k_{-5}) p^3 + \\ &&
   k_{1} k_{4} k_{-4}^3 k_{-5} (k_{5} p^2 + k_{-5}) p^3 + 4 c k_{4} k_{-1} k_{-4}^3 k_{-5}^2 p^2 + 
   4 c k_{1} k_{4} k_{-4}^2 k_{-5} (k_{5} p^2 + k_{-4}) (k_{5} p^2 + k_{-5}) p^2 - \\ &&
   2 c k_{2} k_{5} k_{-4}^2 (k_{4} k_{5}^2 p^6 + 2 k_{4} (k_{4} + k_{5}) k_{-5} p^4 + 
      k_{4} k_{-5} (4 k_{-4} + k_{-5}) p^2 + 2 k_{-4}^2 k_{-5}) S p^2 + \\ && 
   k_{1} k_{-4}^4 (k_{5} p^2 + k_{-5})^2 p + 
   k_{-1} k_{-4}^3 k_{-5} (k_{4} p^2 + k_{-4}) (k_{5} p^2 + k_{-5}) p + \\ &&
   2 k_{-3} k_{-4}^2 (k_{5} p^2 + k_{-5}) (k_{4} k_{5}^2 p^6 + k_{4} (k_{4} + 2 k_{5}) k_{-5} p^4 + 
      k_{4} k_{-5} (2 k_{-4} + k_{-5}) p^2 + k_{-4}^2 k_{-5}) p - \\ &&
   2 k_{3} k_{-4}^2 (k_{5} p^2 + k_{-5}) (k_{4} k_{5}^2 p^6 + k_{4} (k_{4} + 2 k_{5}) k_{-5} p^4 + 
      k_{4} k_{-5} (2 k_{-4} + k_{-5}) p^2 + k_{-4}^2 k_{-5}))
\end{eqnarray*}
\begin{eqnarray*}
\delta_{0}(p) &=& -k_{5} p (-c k_{4} k_{5}^2 k_{-1} k_{-4}^2 k_{-5} p^7 - 2 c k_{4} k_{5} k_{-1} k_{-4}^2 k_{-5}^2 p^5 - \\ &&
   c k_{4} k_{5} k_{-1} k_{-4}^3 k_{-5} p^5 - c k_{4} k_{-1} k_{-4}^2 k_{-5}^3 p^3 - 
   c k_{4} k_{-1} k_{-4}^3 k_{-5}^2 p^3 + \\ &&
   2 c^2 k_{-4}^2 k_{-5} (k_{2} k_{5} k_{-4} (k_{4} p^2 + k_{-4}) S - 
      2 k_{4} k_{-1} k_{-5} (k_{5} p^2 + k_{-4})) p^2 - \\ &&
   c k_{-1} k_{-4}^4 k_{-5} (k_{5} p^2 + k_{-5}) p - 
   2 c k_{-3} k_{-4}^3 k_{-5} (k_{4} p^2 + k_{-4}) (k_{5} p^2 + k_{-5}) p + \\ &&
   2 c k_{3} k_{-4}^3 k_{-5} (k_{4} p^2 + k_{-4}) (k_{5} p^2 + k_{-5})).\\
\end{eqnarray*}
The fast fiber update, which we do not write down, will be an $\mathcal{O}(\eps)$ correction to the first $4\times 2$ block of $A^{(0)}$, in analogy to the previous examples. \\

\begin{remark} The tensorial nature of the update step only becomes apparent in systems with dimension higher than two. For example, to compute $\Lambda_{fs}^{(0)}$, we must compute the tensor $DA_{s}^{(0)}$. Such computations quickly become unwieldy for high-dimensional systems after the first or second CSP iterate.\\
\end{remark}

\begin{remark}
 Symbolic computations for these examples were performed using Mathematica Version 11.1.1.0 \cite{mathematica}. 
\end{remark}

%\newpage
\section{Conclusion}\label{sec:conc}

We have given a detailed treatment of the CSP method for nonstandard slow-fast systems \eqref{eq:master}, when the vector field admits a geometrically intuitive factorization in the leading-order term, and demonstrated the method explicitly for nontrivial examples. There appear to be several other natural connections between the CSP algorithm and this factorization, which we comment on further.

\subsection{CSP as projection} \label{sec:proj}

%\begin{figure}[!]
%  \centering
%    \includegraphics[width=0.8\textwidth,height=0.55\textwidth]{projperp3}
%      \caption{A sketch of the projector operator $\Pi^{\bot}$ defined in \eqref{eq:projperp}. Linear subspaces $T_p S$ and $N_p$, solid lines, are depicted together with their orthogonal complements in dashed lines. }
%      \label{fig:projection2}
%\end{figure}

There is a tantalizing connection between the CSP two-step update and the projection operators \eqref{eq:projS}--\eqref{eq:projN} which were naturally induced by the normally hyperbolic splitting. Consider the CSP fiber of order 1, $\mathcal{L}^{(1)} = \text{Col}(A_f^{(1)})$, defined by \eqref{eq:csp2}. We have
\begin{eqnarray}
A_f^{(1)} &=& A_f^{(0)}(I-\tilde{U}^{(0)}\tilde{L}^{(0)}) + A_s^{(0)} \tilde{L}^{(0)} \nn\\ 
&=& N + A_s^{(0)} \tilde{L}^{(0)} + \mathcal{O}(\eps^2), \label{eq:firstfiber}
\end{eqnarray}

where the term $A_s^{(0)} \tilde{L}^{(0)}$ can be written out explicitly by using the formulas \eqref{eq:lambdaff} and \eqref{eq:lambdafs} for the block components of $\Lambda$ in the update formula \eqref{eq:auxblock}:
\begin{eqnarray}
A^{(0)}_s \tilde{L}^{(0)} &=& (Df^{\top})^{\bot}  ((N^{\bot})^{\top}(Df^{\top})^{\bot})^{-1} (N^{\bot})^{\top} (-DN Nf + \eps (DGN - DN G)) (\Lambda_{11}^{(0)})^{-1} \nn \\
&=& \Pi^{\bot}(-DN Nf + \eps [N,G]) (\Lambda_{ff}^{(0)})^{-1},\label{eq:a1update}
\end{eqnarray}

where 
\begin{eqnarray}
\Pi^{\bot} &:=& (Df^{\top})^{\bot}  ((N^{\bot})^{\top}(Df^{\top})^{\bot})^{-1} (N^{\bot})^{\top}.  \label{eq:projperp}\nn
\end{eqnarray}

A careful comparison to the matrix representation of oblique projection (Definition \ref{def:proj}) reveals that in fact $\Pi^{\bot} = \Pi^{S}$ to leading order, when we evaluate this expression on the set $\mathcal{K}^{(0)}$. Since $f = O(\eps)$ on $\mathcal{K}^{(0)}$ as shown in \eqref{eq:k0}, the equation \eqref{eq:firstfiber} becomes
\begin{eqnarray*}
A^{(1)}_{f,0} &=& N_{0} + \eps \Pi^S_0 \xi_0(N_0,f_0,G_0), \nn
\end{eqnarray*}

where the subscript 0 reminds the reader that these quantities are computed on $\mathcal{K}^{(0)}$, and the vector $\xi_0$ is rigorously defined by the two-step CSP rule \eqref{eq:csp2}. The CSP fiber $\mathcal{L}^{(1)}$ is defined on its corresponding CSP manifold $\mathcal{K}^{(1)}$ by Def. \eqref{eq:cspfibers}, but it is known that the fiber approximation is still $\mathcal{O}(\eps^j)$-accurate if $A_f^{(j)}$ is evaluated on $\mathcal{K}^{(j-1)}$ (see  Sec. 3.5 in \cite{kkz2004b}).  

In the examples \ref{sec:circle} and \ref{sec:modcircle}, $ \Pi^{S}$ is an orthogonal projection onto $T_zS$, and furthermore the manifold $S = S_{\eps}$ is given by $f = 0$ exactly. The computation \eqref{eq:a1update} shows that the fast fiber update is trivial if $[N,G]\in \text{ker~}\Pi^S$.  This is indeed the case in the first example \ref{sec:circle} but not in the modification \ref{sec:modcircle}. 
%is the oblique projection onto $T_pS^{\bot}$ parallel to the orthogonal complement of the linear fast fiber $N_p^{\bot}$  (see Fig. \ref{fig:projection2}).

These results demonstrate the close relationship between the CSP update step and the projectors which define the fast and slow subsystems of \eqref{eq:master}. Deeper connections, including a recasting of the CSP method as a projective scheme in which the projections are themselves iteratively updated, remain to be further clarified.\\

\subsection{CSP as a factorization algorithm for nonstandard vector fields} \label{sec:factor}
 As discussed in Sec. \ref{sec:nongspt}, the leading-order factorization gives a compact description of the geometry of the system near the critical manifold $S$ when $\eps = 0$, and gives approximate information about the dynamics on the nearby slow manifold $S_{\eps}$ when $\eps > 0$. Consider the following conjecture.\\
 
 \begin{conj}
 ~ Vector fields satisfying Assumptions \ref{ass1}--\ref{assu1b} may be re-factorized as follows:
 \begin{eqnarray*}
H(z,\eps) &=& N_{\eps}(z,\eps) f_{\eps}(z,\eps) + \eps \tilde{G}(z,\eps). \label{eq:newfact}
\end{eqnarray*}
 
 The new objects $N_{\eps}(z,\eps)$, $ f_{\eps}(z,\eps)$, and $\tilde{G}(z,\eps)$ satisfy the following properties:
 
 \begin{itemize}
  \item (Slow manifold as a level set) $S_{\eps} = \{z \in \mathbb{R}^n: f_{\eps}(z,\eps) = 0\}$.
 \item (Basis of fast fibers) The columns of $N_{\eps}(p,\eps)$ span the linear fast fibers $\mathcal{N}_{\eps}(p)$ at basepoints $p \in S_{\eps}$.
 \item (Invariance) $\tilde{G}(p,\eps) \in T_p S_{\eps}$ for $p \in S_{\eps}$.
 \end{itemize}
 ~\\
 \end{conj}
The assertion is that a vector field factorization with respect to the critical manifold $S$ when $ \eps = 0$ will induce a new factorization with respect to the slow manifold $S_{\eps}$ when $ \eps > 0$. 

We describe a possible direction to answer this conjecture in the affirmative in the special case where $S_{\eps} = S$ (as in the examples \ref{sec:circle} and \ref{sec:modcircle}). The idea is to use the CSP method to obtain iterative refinements of the objects $N_{\eps}(z,\eps)$, $ f_{\eps}(z,\eps)$, and $\tilde{G}(z,\eps)$. Writing out the first few terms in the asymptotic series of the CSP fast fiber $\mathcal{L}^{(1)}$  (see \eqref{eq:cspfibers} and \eqref{eq:csp2}), we have  $A^{(1)}_f(z,\eps) = N(z) + \eps A^{(1)}_{f,1}(z) + \mathcal{O}(\eps^2)$.

The vector field can be re-factorized as follows:
\begin{eqnarray*}
H(z,\eps) &=& N(z)f(z) + \eps G(z,\eps)\\
&=& N(z) f(z) +  \eps  A^{(1)}_{f,1}(z)  f(z) -  \eps   A^{(1)}_{f,1}(z) f(z) +  \eps G(z,\eps)\\
&=&  (N(z) + \eps A^{(1)}_{f,1}(z)) f(z) + \eps G^{(1)}(z,\eps),\\
&=& N^{(1)}(z,\eps)  f^{(1)}(z) +  \eps G^{(1)}(z,\eps),
\end{eqnarray*}
where
\begin{eqnarray*}
f^{(1)}(z) &:=& f(z)\\
N^{(1)}(z,\eps) &:=& N(z) + \eps A^{(1)}_{f,1}(z)\\
G^{(1)}(z,\eps) &:=& G(z,\eps) - A^{(1)}_{f,1}(z) f(z).
\end{eqnarray*}
In comparison to the original factorization $H = Nf + \eps G$, this new factorization has $\eps$-dependence in the `leading-order' term $N^{(1)}f^{(1)}$, and a modified `remainder' term $G^{(1)}$. But observe that $S = S_{\eps} = \{f^{(1)} = 0\}$. Furthermore, the column space of the prefactor matrix $N^{(1)}$ is now an $\mathcal{O}(\eps)$-approximation of the fast fibers. Finally, $G^{(1)}(z,\eps)|_S = G(z,\eps)  \in T_p S$ since $S$ is invariant by assumption.

We can extend this idea in the obvious way to refine the factorization to some arbitrary order $j$: 
\begin{eqnarray*}
H(z,\eps) &=& N^{(j)}(z,\eps) f^{(j)}(z) + \eps G^{(j)}(z,\eps),\\
f^{(j)}(z) &:=& f(z)\\
N^{(j)}(z,\eps) &:=& N(z) + \sum_{i=1}^j \eps^i A^{(j)}_{f,i}(z)\\
G^{(j)}(z,\eps) &:=& G(z,\eps) - \sum_{i=1}^j \eps^{i-1}A^{(j)}_{f,i}(z) f(z).
\end{eqnarray*}
We remind the reader that the triviality of the level set update $f^{(j)} = f$ is a consequence of the assumption that $S = S_{\eps}$.  The general case where the CSP manifolds also update nontrivially is a topic of further study. Here, more care must be taken to write down the factorizations since the fiber basepoints are not fixed. Furthermore, nontrivial cross-terms of the form $N^{(k)}f^{(l)}$, where $k+l > j$, begin to appear at a given iterate $j$. \\

\subsection{CSP for slow manifolds of saddle-type}

Normally hyperbolic slow manifolds of saddle-type arise in traveling-wave profiles of the FitzHugh-Nagumo model \cite{champneys2007}, in the Hodgkin-Huxley model \cite{hasan2018}, and in models of cardiac pacemakers \cite{krogh2004}. They are very challenging to numerically approximate even in the standard case, due to exponential instabilities in both forward and backward time \cite{GuckKuehn}. It is therefore of interest to relax the assumption that the critical manifold be attracting for the CSP method.\\

~\\

%\newpage

\begin{appendix}
\section{The formalism underlying the CSP method}\label{sec:csp1}

The purpose of this section is to collect a few results given across the three very detailed papers \cite{kkz2004a,kkz2004b,kkz2015}. \\
%
%{\it Changing basis in the variational equation.} Consider an arbitrary smooth vector field $H(z)$, to which we append the variational equation to produce the dynamical system on the tangent bundle:
%
%\begin{eqnarray}
%z' &=& H(z)\\
%H' &=& DH(z) H(z).\nonumber
%\end{eqnarray}
%
%Under some arbitrarily chosen splitting $z = (z_1,z_2)$, where $z_1 \in \mathbb{R}^{d}$ and $z_2 \in \mathbb{R}^{n-d}$ (for some $1 \leq d < n$), we can write the variational equation in block component form:
%
%\begin{eqnarray}
%\begin{pmatrix}
%H_1' \\ H_2'
%\end{pmatrix} &=& 
%\begin{pmatrix}
%D_{z_1} H_1 & D_{z_1} H_2 \\
%D_{z_2} H_1 & D_{z_2} H_2
%\end{pmatrix} \begin{pmatrix}
%H_1 \\ H_2
%\end{pmatrix}. \label{variational2}
%\end{eqnarray}
%
%We study how \eqref{variational} transforms under coordinate changes. Following \cite{kkz2015}, we let $A(z)$ be a smoothly point-dependent $n\times n$ matrix so that the columns of $A(z)$ form a local basis of $T_z \mathbb{R}^n$, and let $B(z)$ be its dual. We define the new vector field $f(z)$
%
%\begin{eqnarray}
%g(z) &=& B(z)H(z),\nonumber
%\end{eqnarray}
%
%to be the {\it vector field $H$ expressed in the basis $A$}.  We can derive several useful identities involving $A, B,H,$ and $g$. Left-multiplying the above by $B(z)$ gives
%\begin{eqnarray}
%H(z) &=& A(z) f(z).\label{eq:ida}
%\end{eqnarray}
%

\subsection{Change of variables and Lie bracket} \label{app:lie} In this section, we justify the statement that the variational equation of the transformed vector field has the structure of a Lie bracket as shown in \eqref{eq:lambdaeq}--\eqref{eq:liebracket}.  Differentiating $BA = I$ with respect to $t$, we obtain the identity

\begin{eqnarray}
BA' &=& - AB'.\label{eq:idb} 
\end{eqnarray}

Furthermore, the chain rule gives

\begin{eqnarray}
A' &=& (DA) z' \label{eq:idc}\\
&=& (DA) H.\nonumber
\end{eqnarray}

Using identities \eqref{eq:idb}-\eqref{eq:idc}, the variational equation of the transformed vector field $f$ becomes:

\begin{eqnarray}
f' &=& (BH)'\nonumber\\
&=& B'H + BH'\nonumber\\
&=& B'Af + B(DH) H\nonumber\\
&=& -BA' f + B(DH) Af\nonumber\\
&=& B(-(DA) H f + (DH) Af)\nonumber\\
&=& B(-(DA) H +  (DH) A) f\nonumber\\
&=& B[A,H] f.
\end{eqnarray}

This calculation appears in \cite{kkz2015}.\\

\subsection{Block-diagonalization of $\Lambda$} \label{app:lambda}
In this section, we justify the characterizations of the invariant manifold and transverse fiber bundle in terms of a cleverly chosen basis, as shown in \eqref{eq:manifold}--\eqref{eq:fiber}. In the presence of an invariant manifold having a transverse fiber bundle which is invariant as a family (i.e. the flow maps fibers into fibers), $H$ can be expressed in a basis which block-diagonalizes the operator $\Lambda$. This was proven in \cite{kkz2004a,kkz2004b} for the case of a slow manifold having a transverse fast fiber bundle, but the identical argument carries over without distinguishing `slow' versus `fast' components. We demonstrate this for the vanishing of the upper-diagonal block. 

Let $\mathcal{M}$ be an invariant manifold with a transverse fiber bundle $\mathcal{N}$; i.e., at a point $z \in \mathcal{M}$, we have the following splitting:
\begin{eqnarray}
T_z \mathbb{R}^n &=& \mathcal{N}_z \oplus T_z \mathcal{M}_z. \nonumber
\end{eqnarray}

We write $A = [A_f ~ A_s]$, where the columns of the $n \times (n-k)$ matrix $A_F(z)$ form a basis of $\mathcal{N}_z$ and the columns of the $n \times k$ matrix $A_s(z)$ form a basis of  $T_z M$. The corresponding dual basis similarly spans the dual splitting: we have $B = \begin{pmatrix} B_{s\bot} \\ B_{f\bot}\end{pmatrix}$. \\

\begin{proposition} ~For the choice of basis $A$ above, we have $\Lambda_{fs} = \mathbb{O}_{n-k,k}$. \end{proposition}
{\it Proof.} We have 

\begin{eqnarray}
\Lambda_{fs} &=& B_{s\bot}[A_s,H]\\
&=& B_{s\bot}[(DH) A_s - (DA_s) H].\nonumber
\end{eqnarray}

By invariance, $H \in TM$ for points in $\mathcal{M}$, and so $B_{s\bot}H = 0$. The directional derivative along $A_s$ must therefore be identically 0 on $\mathcal{M}$ as well. This gives us the following identity:

\begin{eqnarray}
D(B_{s\bot}H) A_s &=& DB_{s\bot}(H,A_s) + B_{s\bot}(DH) A_s= 0. \label{eq:id1}
\end{eqnarray}

Here we clarify that $DB_{s\bot}$ is a symmetric bilinear form, and this (matrix) identity is to be understood as taking the standard directional derivative of each column of $A_s$ in turn and concatenating the result into a matrix.

Similarly, we have the trivial identity $B_{s\bot} A_s = 0$ on $\mathcal{M}$, coming from the dual basis criterion. Differentiating with respect to $t$ and using the chain rule, we obtain

\begin{eqnarray}
\frac{d}{dt} (B_{s\bot} A_s) &=& D(B_{s\bot} A_s) H \label{eq:id2}\\
&=& D(B_{s\bot})(A_s,H) + B_s^{\bot} (DA_s) H = 0.\nonumber
\end{eqnarray}

Subtracting identities \eqref{eq:id1} from \eqref{eq:id2} and using the symmetry of the bilinear form, we obtain $\Lambda_{fs} = 0$. $\Box$

The argument for $\Lambda_{sf} = \mathbb{O}_{n-k,k}$ is similar if slightly more involved, using the invariance of the fiber bundle to obtain an identity as above \cite{kkz2004b}. \\

\end{appendix}

\end{document}